\documentclass[a4paper]{article}

\usepackage{latexsym}
\usepackage{amssymb}
\usepackage{graphicx}

\newcommand{\Doppio}[1]{\mathbb{#1}}

\newcommand{\Tens}[1]{\mathbf{#1}}

\newcommand{\Reali}{\Doppio{R}}

\newcommand{\Euc}{\Doppio{E}_2}
\newcommand{\Min}{\Doppio{M}_2}

\newcommand{\BE}{\begin{equation}}
\newcommand{\EE}{\end{equation}}
\newcommand{\BAS}{\begin{eqnarray*}}
\newcommand{\EAS}{\end{eqnarray*}}
\newcommand{\BA}{\begin{eqnarray}}
\newcommand{\EA}{\end{eqnarray}}

\newcommand{\D}{\mathrm{d}}

\newcommand{\rank}{\mathrm{rank}}
\newcommand{\Mezzo}{\frac{1}{2}}
\newcommand{\Uno}{\mathrm{1\kern-3pt{I}}}

\newtheorem{Prop}{Proposition}

\newtheorem{Lemma}[Prop]{Lemma}
\newtheorem{Cor}[Prop]{Corollary}

\newenvironment{Dim}
{\par\medskip\upshape\noindent\textbf{Proof}:}
{\hfill$\Box$\bigskip}

\newcommand{\Keu}{\mathcal{K}(\Euc)}
\newcommand{\Kmin}{\mathcal{K}(\Min)}

\title{Geometrical classification of Killing tensors on bidimensional flat manifolds}
\author{C. Chanu\footnote{Dipartimento di Matematica, Universit\`a di Torino, Italy},
L. Degiovanni$^*$,
R.G. McLenaghan\footnote{Department of Applied Mathematics,University of Waterloo, Ontario, Canada}}
\date{}

\begin{document}

\maketitle

\begin{abstract}
\noindent Valence two Killing tensors in the Euclidean and Minkowski planes are classified under the action of 
the group which preserves the type of the corresponding Killing web.  The classification is based on an analysis of 
the system of determining partial differential equations for the group invariants and is entirely algebraic. The approach allows one to classify both characteristic and non-characteristic Killing tensors.
\end{abstract}

\section{Introduction and basic properties}
\subsection{Killing tensors and separable webs}
A Killing tensor (KT) on a pseudo-Riemannian space $(M, \Tens{g})$ is a tensor $\Tens{K}$ of type $(0,k)$ which satisfies the equation
\BE\label{lab_1}
\nabla_{(j}K_{i_1 \ldots i_k)} =0\,,
\EE
where $\nabla$ denotes the covariant derivative defined by the Levi-Civita connection of the pseudo-Riemannian metric $\Tens{g}$ and where the parentheses signify symmetrization of the enclosed indices.  It was shown by Eisenhart \cite{ref_1} that such tensors arise naturally from first integrals of the geodesic flow on $(M, \Tens{g})$ in the form
$$
I=K_{i_{1}...i_{k}}\frac{\D q^{i_1}}{\D s}\ldots \frac{\D q ^{i_k}}{\D s}\,.
$$
The function $I$, defined on the the tangent bundle $TM$, is a first integral of the geodesic equations (i.e.\ it is constant along each geodesic) if and only if the Killing tensor equation (\ref{lab_1}) holds.  Killing tensors may also be characterized in contravariant form by means of the following function defined on the cotangent bundle $T^*M$:
$$
I^\ast=K^{i_{1}...i_k}p_{i_1}...p_{i_k}
$$
where $(q^i,p_i)$ denote canonical coordinates on $T^*M$.  Condition (\ref{lab_1}) is then equivalent to
$$
\{I^\ast,H\}=0\,,
$$
where  $\{\cdot,\cdot\}$ denotes the Poisson bracket and
$$
H=\Mezzo\,g^{ij}p_ip_j\,,
$$
the geodesic Hamiltonian.

The set of all Killing tensors of valence $k$, on an $n$-dimensional manifold $M$, is a real vector space which we denote by $\mathcal{K}^k(M)$.  Its dimension $d$  satisfies the Delong-Takeuchi-Thompson inequality \cite{ref_2,ref_3,ref_4}
$$
d\leq \frac{1}{n} {n+k \choose k+1}{n+k-1 \choose k}\,.
$$
Equality is achieved for manifolds of constant curvature.  Moreover, in this case the Killing tensors of valence $k$ are sums of symmetrized products of the Killing vectors of the manifold.  In manifolds with isometry groups of less than the maximal dimension there may exist Killing tensors which are not expressible in this way.  For example, such a situation occurs in the Kerr space-time \cite{ref_5}.
 
Killing tensors of type $(0,2)$ which we call Killing 2-tensors, are particularly important.  Indeed, if the eigenvalues of a Killing 2-tensor are real and simple and the eigenvectors are normal (orthogonally integrable), then the Killing tensor defines an orthogonally separable web on $M$, that is $n$ foliations of mutually orthogonal $(n-1)$-dimensional hypersurfaces.   To the separable web are associated systems of coordinates with respect to which the Hamilton-Jacobi equation for the geodesic flow is solvable by separation of variables (see Benenti \cite{ref_6}). Killing 2-tensors with the above properties are called characteristic Killing tensors (CKT).

It is well known that in the Euclidean plane there exist four types of orthogonally separable webs (see, for example, Miller \cite{ref_7}).  Nevertheless, it is not a trivial task to determine which type of web is defined by a given characteristic Killing tensor.  The converse problem of characterizing the Killing tensors which define the same separable web is also challenging.  This problem becomes even more difficult in dimension greater than two where the preliminary problem of identifying the characteristic Killing tensors is itself a daunting task.  It is thus clear that finding an effective method of classifying Killing tensors would be very useful indeed.

The classification of separable coordinates in two- and three-dimensional 
Euclidean space by Killing tensors dates back to the work of Eisenhart~\cite{ref_1}.  A 
similar classification for two- and three-dimensional Minkowski space was 
undertaken by Kalnins~\cite{Kalnins}, who classified the symmetric second-order 
differential operators that commute with the wave operator and solved the 
Eisenhart integrability conditions~\cite{ref_1} to obtain the metric in the two-dimensional 
case. 
A classification of KT's in the Euclidean and Minkowski planes based on an analysis of their singular sets (i.e.\ the points where the eigenvalues of the Killing tensors are not real and simple) is given by Benenti and Rastelli \cite{ref_8} and Rastelli \cite{ref_9}.  
Recently remarkable progress in the classification problem was achieved 
by McLenaghan, Smirnov, Horwood, The and Yue by means of the Invariant Theory 
of Killing Tensors on spaces of constant curvature \cite{ref_12,ref_10,McL,SmYue}.  In this 
theory Killing tensors are classified modulo the group which consists of the 
transformations on $\mathcal{K}^2(M)$ induced by the isometries of the underlying pseudo-
Riemannian manifold $(M,\Tens{g})$ and the transformation which maps any Killing tensor 
$\Tens{K}$ into $\Tens{K} + b\Tens{g}$, where $b$ is any real number.  More specifically to any isometry 
$\phi$ on $M$ is associated the transformation on $\Tens{K}\mapsto\widehat{\Tens{K}}$ on $\mathcal{K}^2(M)$ defined with 
respect to a system of local coordinates $(q^{i})$ by
\BE\label{Phi_hat}
\widehat{K}^{ij}(q)=J^i_k(\Phi^{-1}(q))J^j_l(\Phi^{-1}(q))K^{kl}(\Phi^{-1}(q))\,,
\EE         
where $J^i_j(q)=\frac{\partial \Phi^i}{\partial q^j}$ is the Jacobian of the transformation $\Phi$.  Two Killing tensors are considered equivalent if one can be obtained from the other in this way.  Clearly all CKT's in the same equivalence class define the same orthogonal web.  
The classification is based on set of 
algebraic invariants of $\mathcal{K}^2(M)$ under the action of the group, from which a 
classification scheme for the type of the separable web can be constructed in 
the cases considered, namely, $\mathbb E_2$, $\mathbb M_2$ and $\mathbb E_3$. 

The approach presented in this paper is related but somewhat different than that developed by McLenaghan et al.  It is based on two observations: (i) the transformations on $\mathcal{K}^2(M)$ induced by the isometries are not the only ones which preserve the type of web defined by a given CKT.  Indeed, any transformation of the form $\Tens{K}\mapsto a\Tens{K}+b\Tens{g}$ also preserves the web (in \cite{ref_10,McL}, $a=1$ 
was assumed); (ii) two webs of the same type are not necessarily isometric.  For example two elliptic-hyperbolic in the Euclidean plane webs with different interfocal distances are of the same type but are not isometric but are rather related by a dilatation transformation.  In the following we do not focus on the transformations of the manifold $M$, but directly on the transformations of $\mathcal{K}^2(M)$ that preserve the type of separable web defined by a given characteristic Killing tensor.  In the Euclidean and Minkowski planes these transformations are well known and generate a Lie group with dimension equal to that of $\mathcal{K}^2(M)$.  This further fact allows the determination of the equivalence classes in a purely algebraic way which is described in the sequel.

There are both advantages and disadvantages to the extension of the group of transformation used in our classification scheme.  On the positive side is the very natural way in which the classes of KT's which define the distinct types of separable webs are obtained.  Restriction
of the transformations of KT's to those that preserve the web has the result that the CKT's which define the same type of web are scattered through many classes.  The method also leaves open the possibility of classifying non-characteristic KT's.  On the negative side, while the isometry group of pseudo-Riemannian manifold is known in many cases, it is not easy to identify the additional transformations of $\mathcal{K}^2(M)$ that preserve the type of a Killing web. This makes it very difficult to extend the method to higher dimensions and to spaces with non-vanishing curvature.

The plan of the paper is as follows: in Subsections 1.2 and 1.3 we outline the necessary theory of Lie transformation groups to be applied later in the paper.  In Section 2 we perform the classification in the Euclidean plane.  The classification in the Minkowski plane is undertaken in Section 3.  Section 4 contains the conclusion.

\subsection{Actions of Lie groups}
Let $A: G\times \mathcal{K} \to \mathcal{K}$ be a linear action of a finite-dimensional Lie
group $G$
on the vector space $\mathcal{K}$. Two points $x,y\in \mathcal{K}$ 
belong to
the same {\em orbit} of the action
if there exists $g\in G$ such that $x=A(g,y)=A_g(y)=g\cdot y$. We call
$\mathcal{O}_x$ the orbit of $A$ containing $x$.
To determine the orbits, we use the infinitesimal generators of the 
action,
i.e.\ the vector fields on $\mathcal{K}$ whose flow coincides with the action 
of  one-parameter subgroups of $G$.
It is well known  that these vector fields form 
a Lie algebra
isomorphic to the Lie algebra of $G$, identified with the set of right-invariant vector fields (see \cite{Olver} for notations).
This means that the distribution spanned by the infinitesimal generators is in fact spanned by just $\dim G$ vector fields (i.e.\ it is \emph{finitely generated}) and it is involutive, but with rank not necessarily constant. We recall that a distribution $\Delta$ is \emph{integrable} if for all $x$ there exists a maximal (connected) integral manifold $S_x$, such that $x\in S_x$ and $\Delta$ is tangent to $S_x$, then $S_x$ is an immersed
submanifold of dimension equal to the rank of $\Delta$ in $x$; the distribution $\Delta$ is \emph{rank-invariant} if the rank of the distribution is constant along the flow of any vector field $X\in\Delta$.

\begin{Prop} \cite{Hermann,Olver}
A distribution $\Delta$ on a manifold is integrable if and only if it is involutive and rank-invariant. Finitely generated involutive distributions are always rank-invariant and hence integrable.
\end{Prop}
The rank-invariance property implies that the distribution $\Delta$ 
is tangent to any of the sets
$$
r_j=\{x\in\mathcal{K}:\rank(\Delta)|_x=j\}\,,
$$
and so if $x\in r_j$ then $S_x\subseteq r_j$, but in general $r_j$ is not a submanifold of $\mathcal{K}$ (not even an immersed one) and is union of several $S_x$. An example is given by the involutive finitely generated distribution on $\Reali^3$ spanned by $\partial_x$ and $(z^2-y^3)\partial_z$, where $r_1$ is not a submanifold.
\begin{Lemma}\label{lemma1}
If $r_j$ is a submanifold of dimension $j$ then for any $x\in r_j$ 
$S_x$ is the connected component of $r_j$ containing $x$.
\end{Lemma}
The proposition follows from the facts that $\Delta$ is tangent to 
$r_j$, $\dim S_x=j$ and $S_x$ is connected.

In our case the distribution $\Delta$ is given by the infinitesimal generator of a Lie group action, then $S_x\subseteq  \mathcal{O}_x$. If $G$ is connected, then its orbits are connected and coincide with the
integral manifolds of $\Delta$. If, instead, $G$ is not connected, then its connected component containing the identity, $G_0$, is a normal 
subgroup.  All the other connected components of $G$ are diffeomorphic to
$G_0$ and coincide with the cosets of $G_0$.
We will denote by $Z$ a set of representatives of the cosets:
$$
G=\bigcup_{g\in Z}g\; G_0.
$$
In all the examples in the following, we can choose $Z$
in such a way that it is a discrete subgroup of $G$.

The orbit $\mathcal {O}_x$ of the action $A$ can be obtained as union of maximal
integral manifolds of $\Delta$ mapped one into the other by the
diffeomorphisms $A_g$ with $g\in Z$
$$
\mathcal{O}_x=\bigcup_{g\in Z}A_g(S_x)=\bigcup_{g\in Z} S_{g\cdot x}.
$$
A consequence of Lemma~\ref{lemma1} is
\begin{Cor}
If $r_j$ is a submanifold of dimension $j$ then for any $x\in r_j$ the 
orbit $\mathcal{O}_x$ is the
union of the connected components of $r_j$ which are images of the one 
containing $x$ through the
action of the elements of $Z$.
\end{Cor}
We conclude by observing that if $\dim G=\dim\mathcal{K}=n$ then we are 
able to determine the orbits
where the distribution $\Delta$ has maximal rank by looking for the 
connected components of
$r_n$ and gluing the ones mapped into the others by the elements of 
$Z$. Moreover  the other
orbits are contained in the sets where the rank of $\Delta$ change and, 
if the condition $\dim
r_j=j$ still holds, they can all be determined in an algebraic way.

\subsection{Sections and connected components}\label{sezione_inco}
The goal of this section is to provide some tools useful to detect the 
components connected by arcs of a subset of $\Reali^n$ (for $n$ big). Actually, in this article with connection we always mean connection by arc, which is equivalent to topological connection in the cases under study.

Let us consider a set $A\subset \Reali^n$ and let
$\{A_i\}_{i\in I}$ be its partition in connected components. 
We consider the natural decomposition $\Reali^n\to \Reali ^m\times \Reali^{n-m}$,
so that any point of $P\in\Reali^n$
can be labeled as
$P=(v,p)$ with $v\in \Reali^m$ and $p\in R^{n-m}$.
In this way we get a partition $\{V^v\}_{v\in\Reali^m}$ of
$\Reali^n$ in parallel hyperplanes of dimension $n-m$, where $V^v=\{P\in \Reali^n: P=(v,p),\; p\in \Reali^{n-m}\}$.
We call $A^v=A \cap V^v$ the section of $A$
determined by the hyperplane $V^v$ and construct its partition in 
connected components
$$
A^v=\bigcup_{\alpha\in I^v} A_\alpha^v\,.
$$
On the family of the connected components of the sections of $A$:
$$
\left\{A_\alpha^v\right\}_{{\alpha\in I^v},\,{v\in \Reali^m}}
$$
we define the relation
$$
A_\alpha^v \sim A_\beta^w \iff \exists!\; i\in I: A_\alpha^v \subseteq
A_i\mbox{ and  }A_\beta^w\subseteq A_i
$$
It's easy to check that the following Lemma holds:
\begin{Lemma}\label{inco}
The relation $\sim$ is an equivalence relation. There is a one-to-one 
correspondence between the
equivalence classes $\left[A_\alpha^v\right]$ and the (arc)-connected 
components $A_i$ of $A$. If there
exists a continuous arc $f:[0,1]\to\Reali^m\times\Reali^{n-m}$ 
such that $f([0,1])\subseteq
A$, $f(0)=(v,p)$, $f(1)=(w,q)$ with $p\in A_\alpha^v$, $q\in A_\beta^w$ 
then $A_\alpha^v\sim
A_\beta^w$.
\end{Lemma}
 From this Lemma it follows that the study of the connected components 
of $A$ can be reduced to the
study of connected components of all sections $A^v$ (of lower 
dimension) under the equivalence
relation.

\section{Killing tensors in the Euclidean plane}
In the Euclidean plane $\Euc$, with Cartesian coordinates $(x,y)$ and  
the standard metric $\Tens{g}$, the general Killing 2-tensor $\Tens{K}$  
has the following contravariant form:
$$
\|K^{ij}\|=\left(\begin{array}{cc}
A+2\,\alpha\,y+\gamma\,y^2 & C-\alpha\,x-\beta\,y-\gamma\,xy \\
C-\alpha\,x-\beta\,y-\gamma\,xy & B+2\,\beta\,x+\gamma x^2
\end{array}\right)\,.
$$
We denote by $\Keu$ the vector space of KT's on the Euclidean plane. On this  
space there exist six kinds of transformation preserving the type of the  
web associated to each KT: three of them correspond to  
isometries, a fourth corresponds to the dilatation of $\Euc$. The  
last two are not associated with any coordinate transformation in the  
plane but act directly on the tensor $\Tens{K}$ and correspond to the  
addition of a multiple of the metric tensor  
($\Tens{K}\mapsto\Tens{K}+\tau\Tens{g}$) and to the multiplication of  
the tensor for a non-vanishing constant  
($\Tens{K}\mapsto\lambda\Tens{K}$). The infinitesimal generators of  
these transformations are easily calculated (see \cite{ref_10} for the  
generators corresponding to isometries and addition of a multiple of the 
metric). With respect to the basis of  
the vector fields on $\Keu$ given by  
$(\partial_A,\partial_B,\partial_C,\partial_\alpha,\partial_\beta,\partial_\gamma)$ the infinitesimal generators are spanned by:\\

Translations
\BAS
V_1 &= & \left(0, -2\beta, \alpha, 0, -\gamma, 0\right) \\
V_2 &=&  \left(-2\alpha, 0, \beta, -\gamma, 0, 0\right)
\EAS

Rotation
$$
V_3 =  \left(-2C, 2C, A-B, \beta, -\alpha, 0\right)
$$

Dilatation of $\Euc$
$$
V_4 =  \left(2A, 2B, 2C, \alpha, \beta, 0\right)
$$

Addition of the metric
$$
V_5 =  \left(1, 1, 0, 0, 0, 0\right)
$$

Scalar multiplication
$$
V_6 =  \left(A, B, C, \alpha, \beta, \gamma\right)
$$

These vector fields form a Lie algebra and therefore generate an  
integrable distribution, denoted by $\Delta_E$. In order to study  
the rank of $\Delta_E$ we gather the components of the $V_i$ in the  
matrix
\BE\label{mat_euc}
M=\left(
\begin{array}{cccccc}
0 & -2\beta & \alpha &  0 &  -\gamma &  0 \\
-2\alpha &  0 &  \beta &  -\gamma &  0 &  0 \\
-2C & 2C &  A-B &  \beta &  -\alpha &  0 \\
2A &  2B &  2C & \alpha &  \beta &  0 \\
1 & 1 &  0 &  0 &  0 &  0 \\
A &  B &  C &  \alpha & \beta &  \gamma
\end{array}
\right)
\EE
with determinant:
$$
\det M =  
-2\gamma\,\left[(\alpha^2-\beta^2-\gamma(A-B))^2+4(\alpha\beta+\gamma  
C)^2\right]\,.
$$
We are led naturally to consider the two surfaces where $\det M=0$
\BE
S_1: \gamma=0 \qquad \mbox{dim }S_1=5
\EE
\BE\label{S2}
S_2: \left\{
\begin{array}{l}
\alpha^2-\beta^2 = \gamma(A-B) \\
\alpha\beta = -\gamma C
\end{array}
\right.
\qquad \mbox{dim }S_2=4\,.
\EE
whose intersection is the vector subspace $\alpha=\beta=\gamma=0$.

The sections of $S_1$, obtained using as parameters $A$, $B$ and $C$,  
are always planes; on the other hand the sections of $S_2$ are curves  
described by the following lemma:
\begin{Lemma}\label{forma_S2}
If the parameters $A$, $B$ and $C$ have the values $C=0$ and $A=B$ then  
the section of $S_2$ is given by the axis $\gamma$, for other values of  
the parameters the section is given by two parabolas contained  
in two orthogonal planes, with  
vertex in the origin and foci on the $\gamma$ axis symmetric with respect to the origin.
\begin{Dim}
Firstly we consider the case $C\neq0$, then the equations (\ref{S2})  
can be transformed into
$$
\left\{
\begin{array}{l}
(\alpha^2-\beta^2)C+\alpha\beta(A-B) = 0 \\
\alpha\beta = -\gamma C
\end{array}
\right.
$$
The first equation can be factorized as  
$C(\alpha-k_+\beta)(\alpha-k_-\beta)$ where
$$
k_\pm=\frac{B-A \pm \sqrt{(A-B)^2+4C^2}}{2C}\,.
$$
Thus the section of $S_2$ is the union of the two parabolas
$$
\left\{
\begin{array}{l}
\alpha=k_+\beta \\[7pt]
\gamma=-\displaystyle{\frac{k_+}{C}}\beta^2
\end{array}
\right.
\bigcup\;
\left\{
\begin{array}{l}
\alpha=k_-\beta \\[7pt]
\gamma=-\displaystyle{\frac{k_-}{C}}\beta^2
\end{array}
\right.
$$
We observe that $k_+k_-=-1$ thus the two parabolas  
are contained in two orthogonal planes. Their foci lie on the $\gamma$ axis with
$$
\gamma=\pm\frac{1}{4}\sqrt{(A-B)^2+4C^2}\,,
$$
being $k_+/C>0$ the first parabola is always downward, while the second one is always upward.
For $C=0$ the second equation in (\ref{S2}) becomes $\alpha\beta=0$.  
Then when $A\neq B$ we have the two parabolas
$$
\left\{
\begin{array}{l}
\alpha=0 \\[5pt]
\gamma=\displaystyle{\frac{\beta^2}{B-A}}
\end{array}
\right.
\bigcup\;
\left\{
\begin{array}{l}
\beta=0 \\[5pt]
\gamma=\displaystyle{\frac{\alpha^2}{A-B}}
\end{array}
\right.
$$
for which the previous considerations on foci hold. Finally when $A=B$ we have $\alpha=\beta=0$ and then the two parabolas degenerate in the $\gamma$ axis.
\end{Dim}
\end{Lemma}

We remark that the functions $\gamma$ and  
$\delta:=(\alpha^2-\beta^2-\gamma(A-B))^2+4(\alpha\beta+\gamma C)^2$, defining  
the surfaces $S_1$ and $S_2$ are the
fundamental invariant of $\mathcal K_2(\mathbb R^2)$  determined by  
McLenaghan et al. \cite{ref_10} under the action of 
the group induced by the isometries and the addition of a multiple of the 
metric.

\begin{Prop}
Outside of the union of the surfaces $S_1$ and $S_2$, the distribution  
$\Delta_E$ has rank $6$ and the space $\Keu-(S_1 \cup S_2)$ is an orbit  
of the action.
\begin{Dim}
The determinant of the matrix (\ref{mat_euc}) is
$$
\det M =  
-2\gamma\,\left[(\alpha^2-\beta^2-\gamma(A-B))^2+4(\alpha\beta+\gamma  
C)^2\right]\,.
$$
Hence, the distribution has maximal rank outside of $S_1 \cup S_2$.  
Since $\Keu-(S_1 \cup S_2)$ has the same dimension of the distribution,  
each connected component is an orbit of the action generated by the  
vector fields $V_i$. The connected components are two: one for  
$\gamma>0$ and the other for $\gamma<0$. However, the two components  
are linked together by the finite transformation that change the sign  
of the KT and so they form a unique orbit with respect to  
the disconnected group generated by the vector fields and this  
transformation.
\end{Dim}
\end{Prop}

\begin{Prop}
On $S_1-S_2$ the rank of the distribution $\Delta_E$ is $5$ and this  
space is an orbit of the action.
\begin{Dim}
In order to determinate the rank of $\Delta_E$ on $S_1$ we set  
$\gamma=0$ in the matrix $M$ and look at its $5\times 5$ minors. This  
task can be easily performed calculating the adjoint matrix of $M$:
$$
\mathrm{adj}(M)|_{\gamma=0}= (\alpha^2+\beta^2)\,
\left(
\begin{array}{cccccc}
0&0&0&0&0&0\\
0&0&0&0&0&0\\
0&0&0&0&0&0\\
0&0&0&0&0&0\\
0&0&0&0&0&0\\
\ast & \ast & 0 & 2(\alpha^2+\beta^2) & \ast & -2(\alpha^2+\beta^2)
\end{array}
\right)\,.
$$
This matrix vanishes identically (and so the rank of $M$ is lesser than  
5) if and only if $\alpha=\beta=0$, that is on $S_1 \cap S_2$. Because  
$S_1$ without its intersection with $S_2$ is connected and it has  
dimension equal to the rank of the distribution on it, then $S_1-S_2$  
is an orbit of the action.
\end{Dim}
\end{Prop}

\begin{Prop}
On $S_2-S_1$ the rank of the distribution $\Delta_E$ is $4$ and this  
space is an orbit of the action.
\begin{Dim}
Assumed $\gamma\neq0$, from the equations (\ref{S2}) we obtain the  
relations
$$
B=A-\frac{\alpha^2-\beta^2}{\gamma}\,, \qquad  
C=-\frac{\alpha\beta}{\gamma}
$$
which substituted in $\mathrm{adj}(M)$ make it identically zero. Then  
on $S_2-S_1$ the rank of the distribution is at most 4, but  the  
$4\times 4$ minor of $M$ obtained by eliminating the second and third  
columns and the third and forth rows is
$$
\left| \begin{array}{cccc}
0 &0 & -\gamma &0\\
-2\alpha & -\gamma & 0 &0\\
1 & 0 & 0 & 0\\
A & \alpha & \beta & \gamma
\end{array}\right| =\gamma^3\neq0
$$
and so outside of $S_1 \cap S_2$ the rank is exactly 4. From  
Lemma~\ref{forma_S2} it follows that for any fixed values of $A-B$ and $C$  
(not both vanishing) the section of $S_2-S_1$ is formed by four  
disjoint parabola's arcs. But it is always possible to find a continuous  
deformation of the parameter $A,B,C$ gluing together the two upward and  
downward arcs, respectively. Indeed with the change in the space of  
parameters $A-B=\rho\cos\theta$, $2C=\rho\sin\theta$ we have that the  
directions of the two planes containing the parabolas depends only on  
$\theta$, while the amplitude of the two parabolas is inversely  
proportional to $\rho$, thus letting $\rho$ go to zero with a fixed  
value of $\theta$ has the effect to glue together the arcs of the two  
parabolas along the $\gamma$ axis. Hence, $S_2-S_1$ has two connected  
components only which can be connected using the change of sign of the  
KT.
\end{Dim}
\end{Prop}

Finally we study the intersection $S_1 \cap S_2$ which is the  
three-dimensional vector space with coordinates $A$, $B$ and $C$. On  
$S_1 \cap S_2$ the only independent vector fields among the $V_i$ are  
$V_3$, $V_5$ and $V_6$, whose components, with respect to  
$(\partial_A,\partial_B,\partial_C)$, form the matrix
\BE\label{mat_euc_rid}
\widetilde{M}=\left(
\begin{array}{ccc}
-2C & 2C & A-B \\
A & B & C \\
1 & 1 & 0
\end{array}
\right)
\EE
Introducing the one-dimensional line
\BE
S_3: \left\{
\begin{array}{l}
\alpha = \beta = \gamma = 0 \\
C = 0 \\
A = B
\end{array}
\right.
\qquad \mbox{dim }S_3=1
\EE
we are able to individuate the last two orbits.

\begin{Prop}
The rank of the distribution $\Delta_E$ on $(S_1 \cap S_2) - S_3$ is  
$3$ and then this space is an orbit of the action.
\begin{Dim}
The determinant of the matrix (\ref{mat_euc_rid}) is $\det  
\widetilde{M} =4C^2 + (A-B)^2$, then it vanishes only on $S_3$. Because  
$(S_1 \cap S_2) - S_3$ is connected it is an orbit.
\end{Dim}
\end{Prop}

\begin{Prop}
The rank of the distribution $\Delta_E$ on $S_3$ is $1$ and then this  
space is an orbit of the action, containing the (non-characteristic)  
tensors of the form $\tau\,\Tens{g}$.
\begin{Dim}
The only independent vector field on $S_3$ is the constant vector  
$V_5$, generated by the addition of a multiple of the metric.
\end{Dim}
\end{Prop}

We remark that the discrete transformation $(A\leftrightarrow B, \alpha\leftrightarrow\beta)$ induced by the discrete isometry of the Euclidean plane ($\{\bar{x}=y, \bar{y}=x\}$) does not allow one to glue together the above found orbits.

In conclusion five orbits of the action of the web preserving group are  
found.
\begin{description}
\item{E1)} The set $\Keu-(S_1 \cup S_2)$, the tensors on this orbit generate  
elliptic-hyperbolic coordinates. A tensor of this type is:
$$
\left(
\begin{array}{cc}
y^2 & 1-xy  \\
1-xy & x^2
\end{array}
\right)\,.
$$
\item{E2)} The set $S_1-S_2$, the tensors on this orbit generate parabolic  
coordinates. Two tensors of this type are:
$$
\left(
\begin{array}{cc}
2y & -x  \\
-x & 0
\end{array}
\right)
\quad \mbox{and}\quad
\left(
\begin{array}{cc}
0 & -y  \\
-y & 2x
\end{array}
\right)\,.
$$
\item{E3)} The set $S_2-S_1$, the tensors on this orbit generate polar coordinates. A tensor of this type is:
$$
\left(
\begin{array}{cc}
y^2 & -xy  \\
-xy & x^2
\end{array}
\right)\,.
$$
\item{E4)} The set $(S_1 \cap S_2)-S_3$, the tensors on this orbit generate Cartesian  
coordinates. Three tensors of this type are:
$$
\left(
\begin{array}{cc}
1 & 0  \\
0 & 0
\end{array}
\right),
\quad
\left(
\begin{array}{cc}
0 & 0  \\
0 & 1
\end{array}
\right),
\quad\mbox{and}\quad
\left(
\begin{array}{cc}
0 & 1  \\
1 & 0
\end{array}
\right)\,.
$$
\item{E5)} The line $S_3$, the tensors on this orbit are multiples of the metric.
\end{description}
This classification coincides with that given by McLenaghan et al. \cite{ref_10} where 
the four types of separable webs in $\mathbb E_2$ are characterized by the vanishing or 
not of the fundamental invariants $\gamma$ and $\delta$.
The orbits are strictly related to the set of singular points discussed by Benenti and Rastelli \cite{ref_8}.
Indeed, the discriminant of the characteristic polynomial of $\Tens{K}$ vanishes on points satisfying
\BE \label{discreu}
\left \{
\begin{array}{l}
\gamma xy+\alpha x+\beta y-C=0 \\
\gamma(y^2-x^2)+2(\alpha y-\beta x)+A-B=0
\end{array}
\right.
\EE
If $\gamma \neq 0$ (i.e.\ outside $S_1$), the equations (\ref{discreu}) describe two hyperbolas both centered in $(-\frac{\beta}{\gamma},-\frac{\alpha}{\gamma})$.
For tensors belonging to $S_2-S_1$ both conics degenerate into two couples of lines through the center
(polar web). Otherwise, they have two points in common (elliptic-hyperbolic web).
For tensor belonging to $S_1$ ($\gamma=0$)  the system (\ref{discreu}) is linear: if $\Tens{K}\in S_1-S_2$, it  represents the intersection of two orthogonal lines (parabolic web); if $\Tens{K}\in (S_1\cap S_2)-S_3$ the system has no solution (Cartesian web), while for tensors belonging to $S_3$ all points are singular. 

\section{Killing tensors in the Minkowski plane}
On the Minkowski plane $\Min$ with pseudo-Cartesian coordinates $(t,x)$  
and metric $\Tens{g}$ with contravariant components
$$
\|g^{ij}\|=\left(
\begin{array}{cc}
1 & 0 \\
0 & -1
\end{array}
\right)
$$
the general Killing tensor $\Tens{K}$ has contravariant components:
$$
\|K^{ij}\|=\left(\begin{array}{cc}
A+2\,\alpha\,x+\gamma\,x^2 & C+\alpha\,t+\beta\,x+\gamma\,tx \\
C+\alpha\,t+\beta\,x+\gamma\,tx & B+2\,\beta\,t+\gamma t^2
\end{array}\right)\,.
$$
We denote by $\Kmin$ the vector space of all the KT's on  
$\Min$. On this space six kinds of transformation are defined, which preserve the type of the web associated to the KT: 
three are induced by the isometries of the Minkowski plane and a fourth by its 
dilatation; the last two do not depend on any transformation of  $\Min$ and are defined directly on $\Kmin$: adding a  
multiple of the metric tensor ($\Tens{K}\mapsto\Tens{K}+\tau\Tens{g}$)  
and multiplying the tensor for a non-vanishing constant  
($\Tens{K}\mapsto\lambda\Tens{K}$). With respect to the basis of the  
vector fields on $\Kmin$ given by  
$(\partial_A,\partial_B,\partial_C,\partial_\alpha,\partial_\beta,
\partial_\gamma)$ the infinitesimal generators are spanned by:\\

Translations:
\BAS
V_1 &= & \left(0, -2\beta, -\alpha, 0, -\gamma, 0\right) \\
V_2 &=&  \left(-2\alpha, 0, -\beta, -\gamma, 0, 0\right)
\EAS

Boost (hyperbolic rotation)
$$
V_3 =  \left(2C, 2C, A+B, \beta, \alpha, 0\right)
$$

Dilatation of $\Min$
$$
V_4 =  \left(2A, 2B, 2C, \alpha, \beta, 0\right)
$$

Addition of the metric
$$
V_5 =  \left(1, -1, 0, 0, 0, 0\right)
$$

Scalar multiplication
$$
V_6 =  \left(A, B, C, \alpha, \beta, \gamma\right)
$$
(see \cite{McL} for the computation of $V_1$, $V_2$, $V_3$, and $V_5$).

Moreover, similar to the Euclidean case, there are the 
following discrete transformations which are analyzed in detail in subsection 
~\ref{sezione_discreta}: the first is the change in sign of the Killing tensor
$$
R_0: K \to -K.
$$                      
The others are
induced from the discrete isometries of $\Min$ $\{  
\bar t=t,\;\bar x=-x\}$
and $\{\bar t=-t,\;\bar x= x\}$, they are
$$
R_1: C\to-C, \alpha\to-\alpha
$$
$$
R_2: C\to-C, \beta\to-\beta
$$
In \cite{Kalnins} and\cite{McL} the transformations used are $R_1$ together with
$$
\widehat{R_2}: A\leftrightarrow B,  \alpha\leftrightarrow\beta\,,
$$
which arises from a change of signature of the metric.
We prefer transformation $R_2$ instead of $\widehat{R_2}$ because it preserves the 
interior (and exterior) of the null cone in $\Min$.

\subsection{Study of the distribution rank}\label{sezione_continua}
The vector fields $V_i$ form a Lie algebra and therefore generate an  
integrable distribution, denoted by $\Delta_M$. In to order to study  
the rank of $\Delta_M$ we gather the components of the $V_i$ in the  
matrix
\BE\label{mat_min}
M=\left(
\begin{array}{cccccc}
0 &  -2\beta & -\alpha & 0 & -\gamma & 0 \\
-2\alpha & 0 & -\beta & -\gamma & 0 & 0 \\
2C & 2C & A+B & \beta & \alpha & 0 \\
2A &  2B &  2C & \alpha &  \beta &  0 \\
1 & -1 &  0 &  0 &  0 &  0 \\
A &  B &  C &  \alpha & \beta &  \gamma
\end{array}
\right)
\EE
with determinant:
$$
\det M =  
2\gamma\,\left[\gamma(A+B-2C)-(\alpha- 
\beta)^2\right]\,\left[\gamma(A+B+2C)-(\alpha+\beta)^2\right]\,.
$$
Thus we consider the two surfaces
\BE
S_1: \gamma=0 \qquad \mbox{dim }S_1=5
\EE
\BE\label{S2m}
S_2:  
\Big[\gamma(A+B-2C)-(\alpha-\beta)^2\Big]\,\Big[\gamma(A+B+2C)- 
(\alpha+\beta)^2\Big]=0 \qquad \mbox{dim }S_2=5.
\EE
We remark that the functions
$$
f_1=\gamma, \quad f_2=
\Big[\gamma(A+B-2C)-(\alpha-\beta)^2\Big]\,\Big[\gamma(A+B+2C)- 
(\alpha+\beta)^2\Big],
$$
coincide with the two fundamental algebraic invariants of $\Kmin$ under the
action of the isometry group augmented by addition of a multiple of the metric
given in \cite{McL}.

The surface $S_2$ is formed by two branches $B_1$ and $B_2$ given,  
respectively, by the equations $\gamma(A+B-2C)=(\alpha-\beta)^2$ and  
$\gamma(A+B+2C)=(\alpha+\beta)^2$.  Nevertheless these two branches are  
mapped one in the other by the transformation $R_1$ and thus it is  
appropriate to consider them as a unique object.  The intersection of  
$B_1$ and $B_2$ is the surface
\BE\label{S3m}
B_1 \cap B_2 = S_3: \left\{
\begin{array}{l}
\gamma(A+B)=\alpha^2+\beta^2 \\
\gamma C = \alpha\beta
\end{array}
\right.
\qquad \mbox{dim }S_3=4
\EE
The intersection of $S_1$ and $S_2$ is described by the equations  
$\gamma=0$ and $\alpha^2=\beta^2$, while $S_1 \cap S_3$ has equations  
$\alpha=\beta=\gamma=0$.

\begin{figure}
\includegraphics[width=5cm]{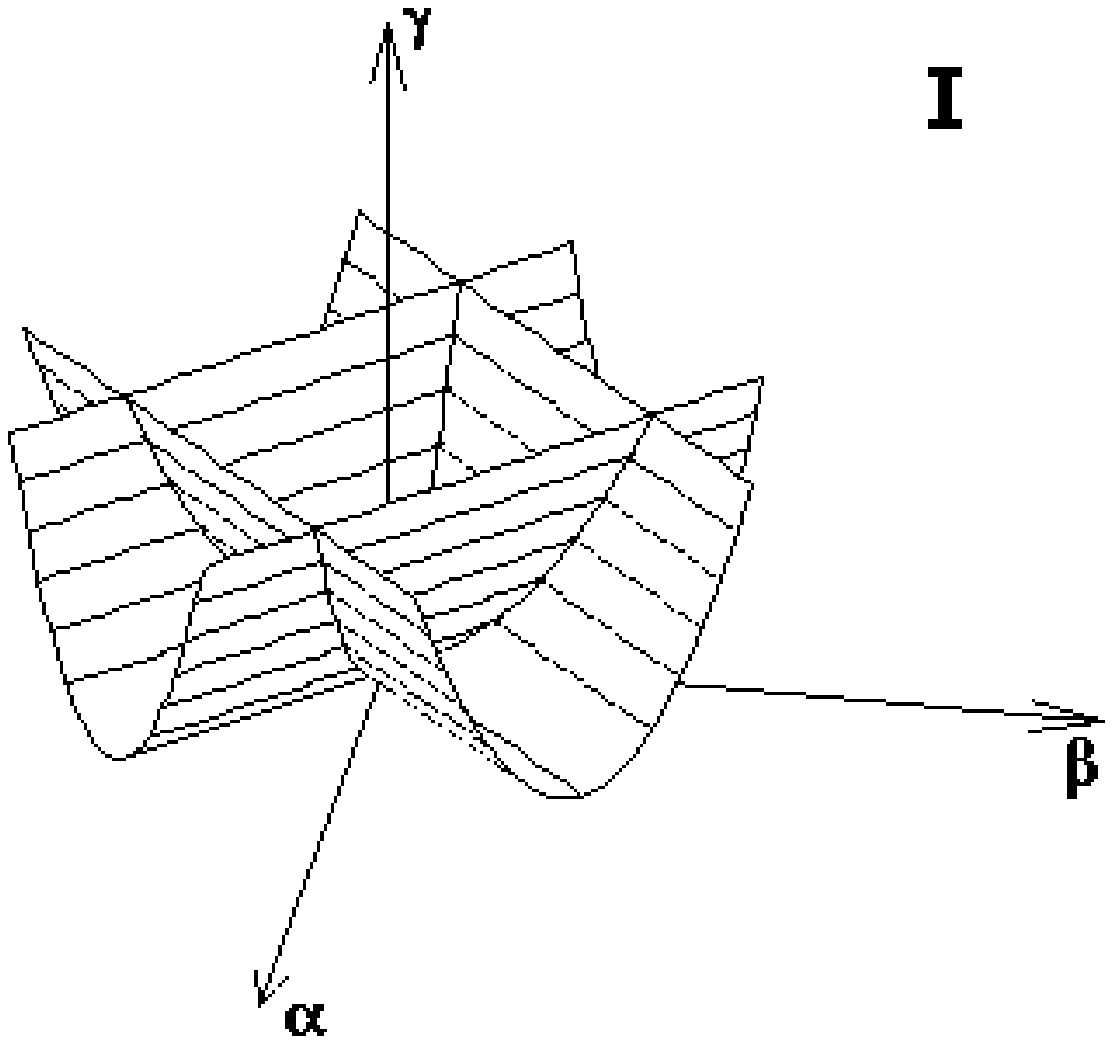}
\includegraphics[width=5cm]{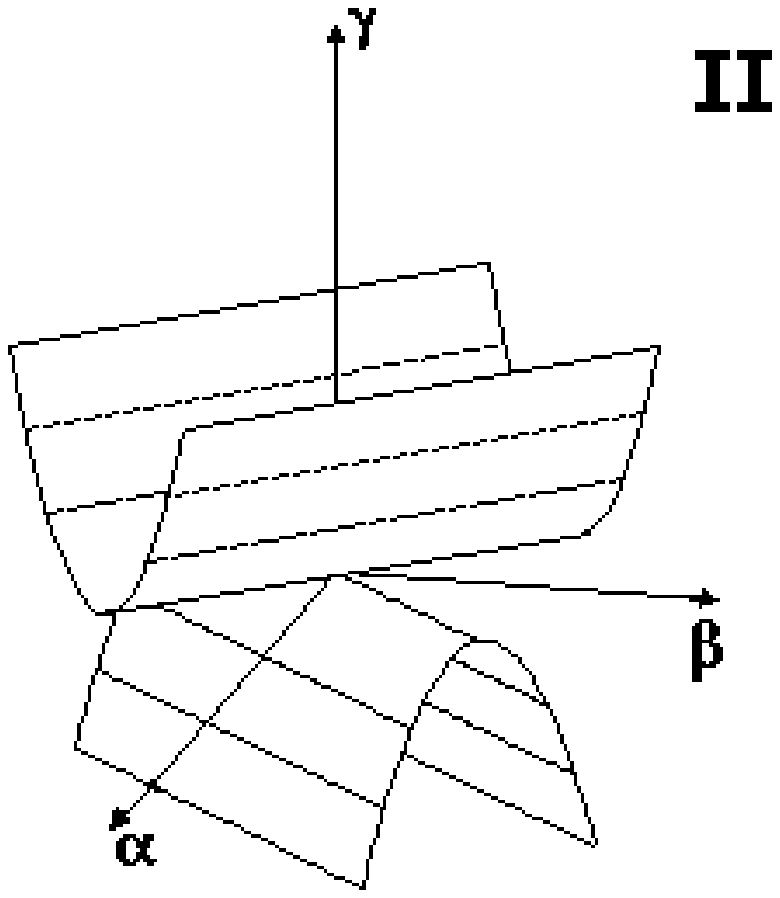}
\includegraphics[width=5cm]{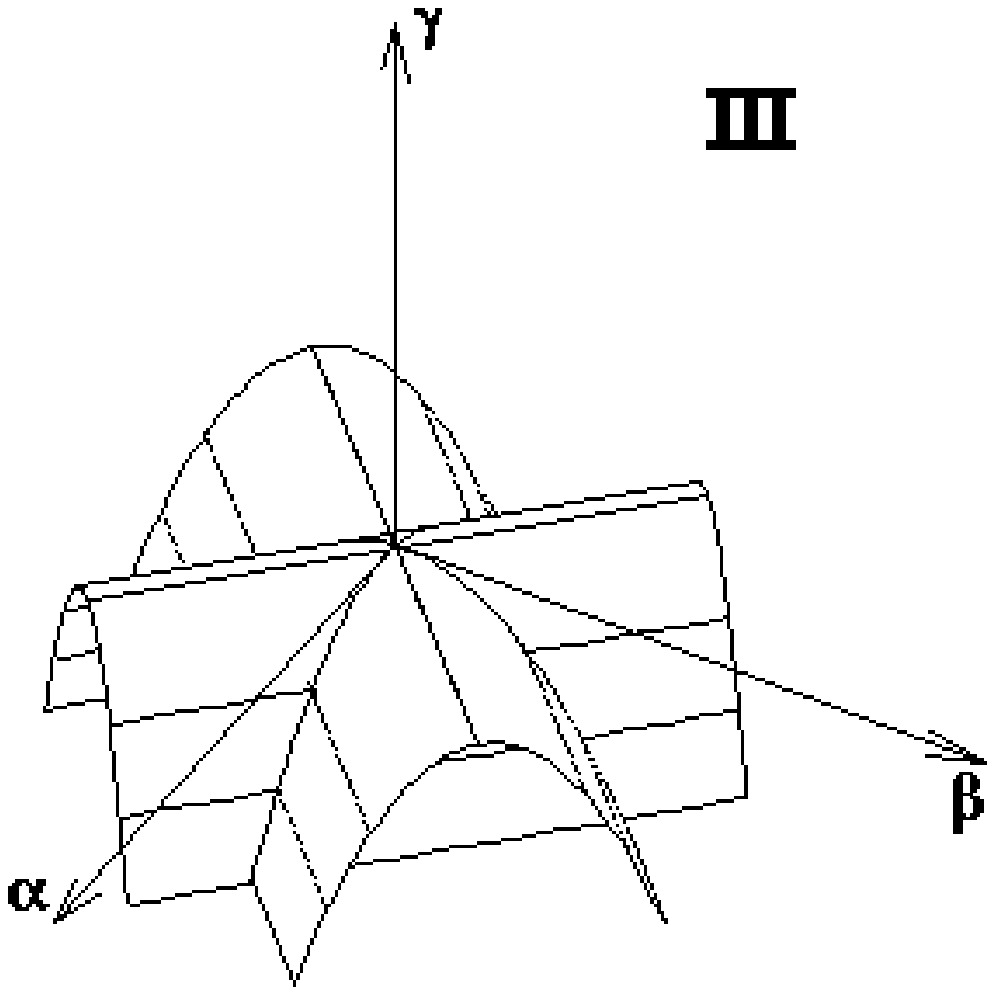}
\includegraphics[width=5cm]{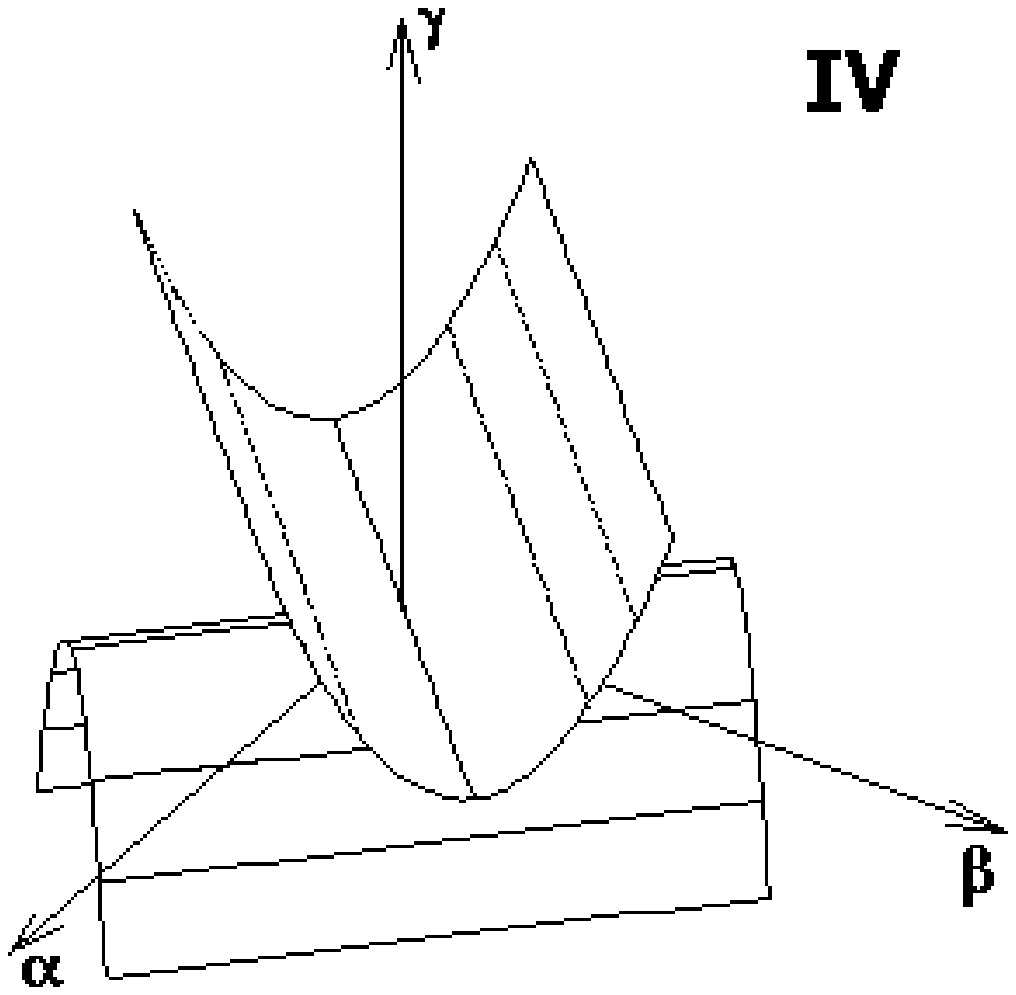}
\includegraphics[width=7cm]{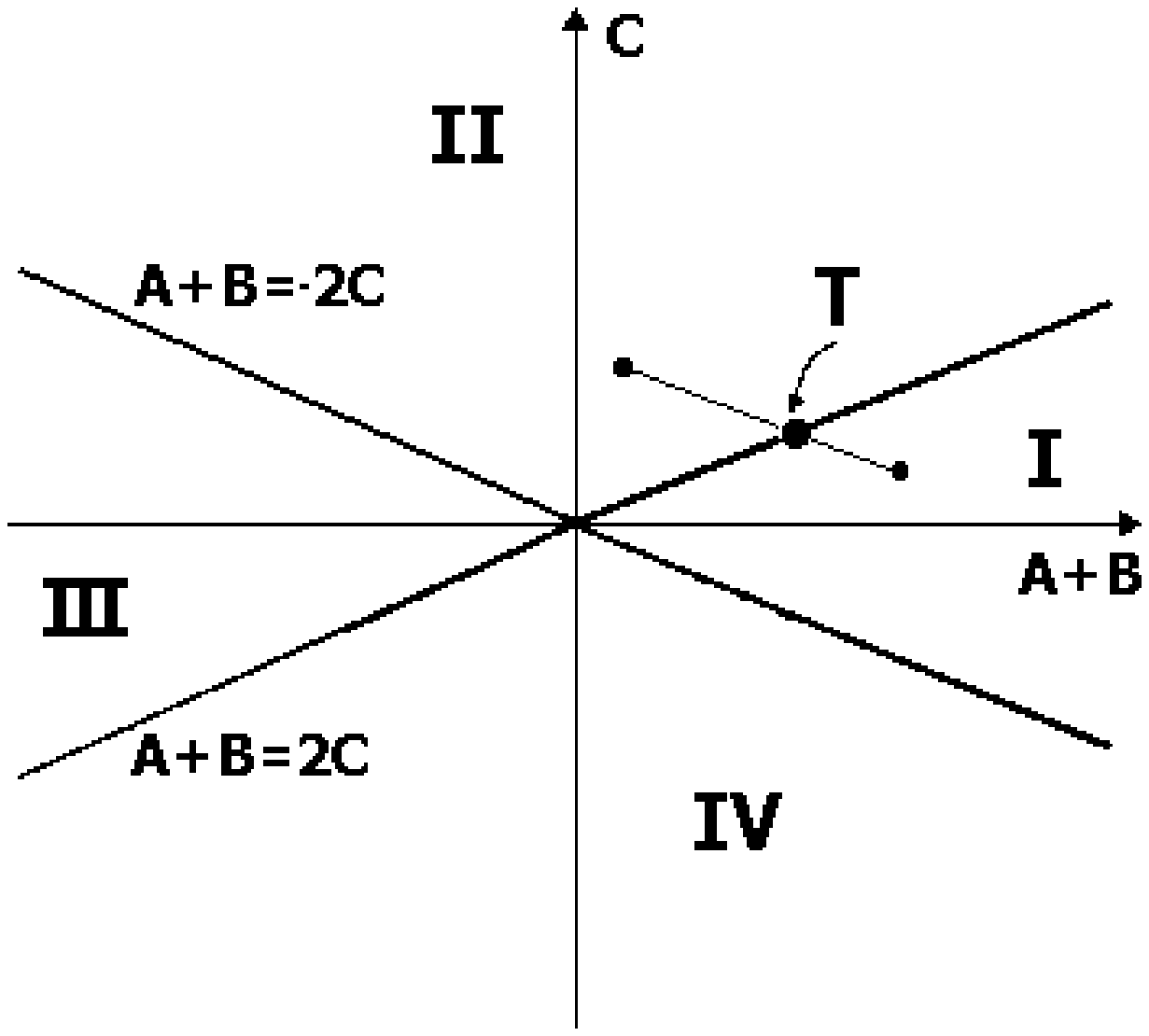}
\includegraphics[width=5cm]{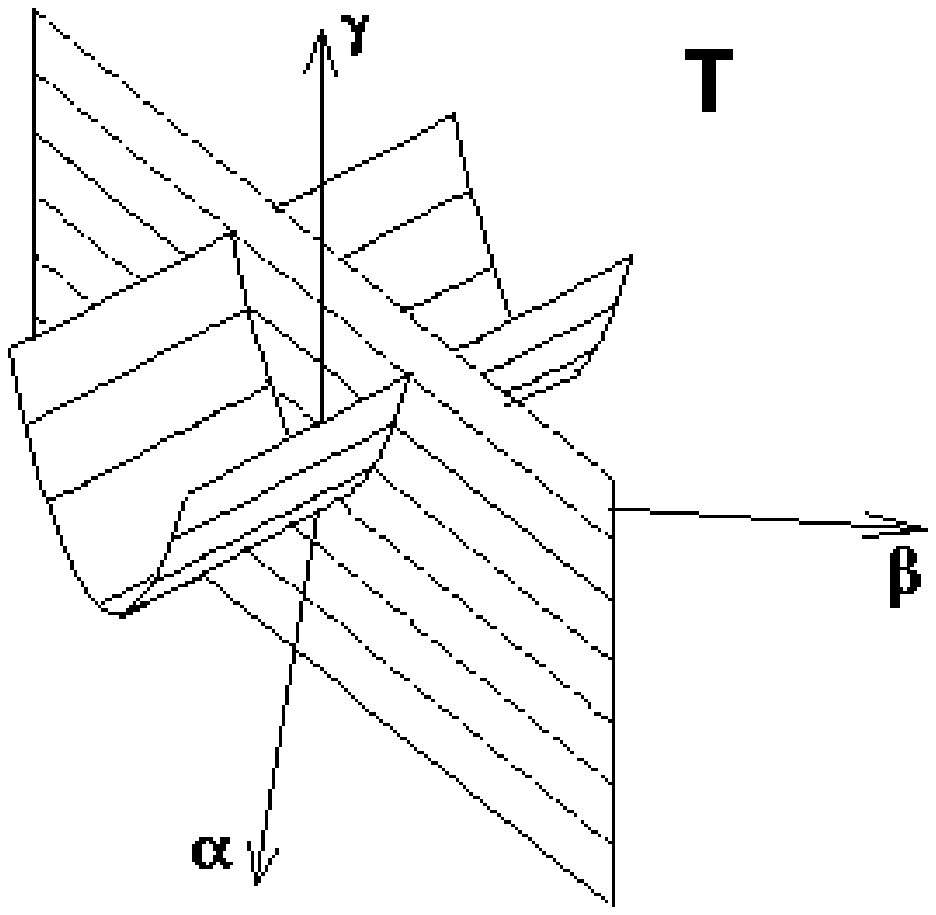}
\caption{The sections of surface $S_2$}\label{pic_S2}
\end{figure}

\begin{figure}
\centering
\includegraphics[width=10cm]{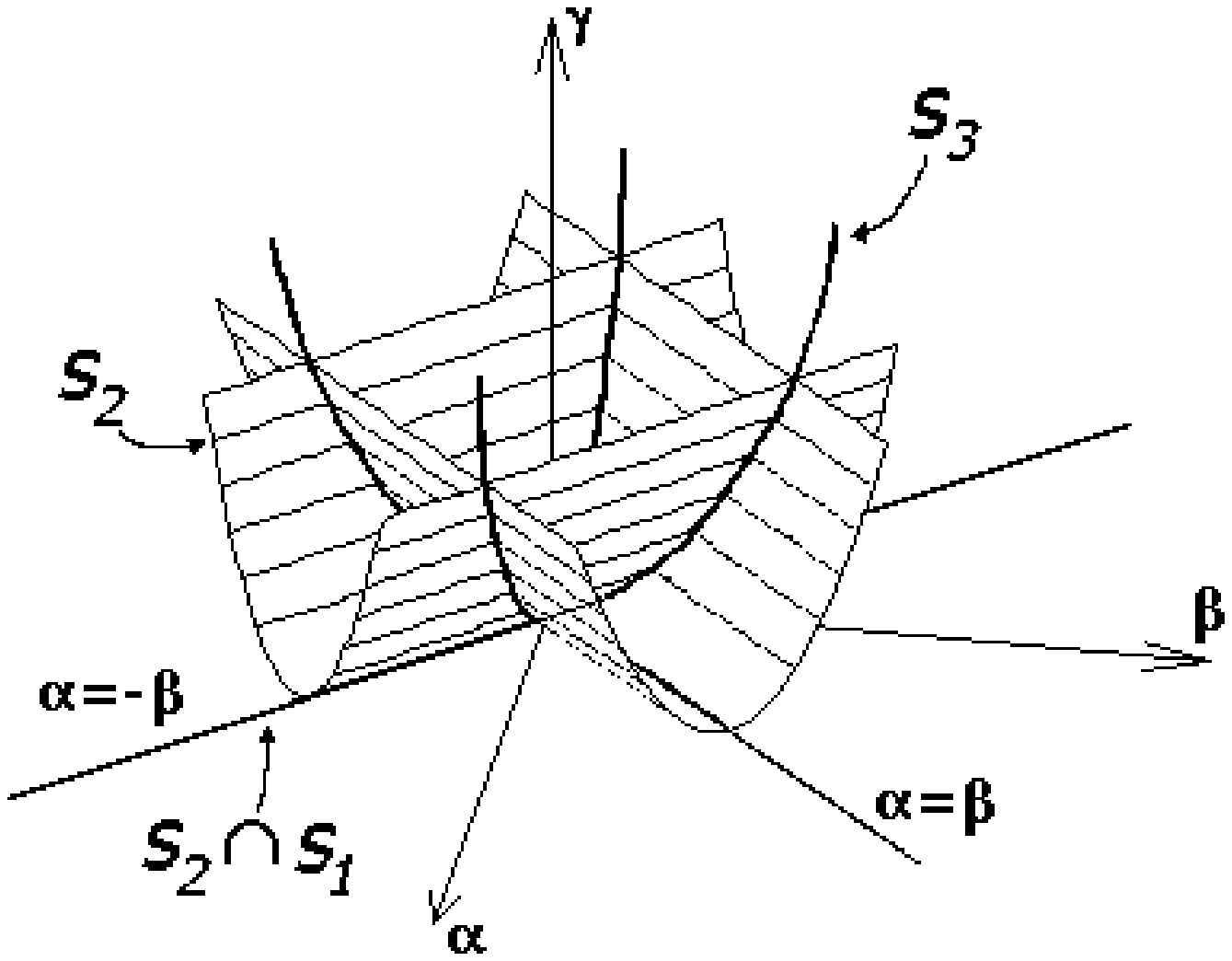}
\caption{The section of $S_1$, $S_2$, $S_3$ in the region  
I}\label{pic_Svari}
\end{figure}

The surfaces in Figure~\ref{pic_S2} represent all the possible  
(generic) sections of $S_2$ in the space $\alpha,\beta,\gamma$. These  
sections can be grouped in four kinds, corresponding to the following  
open sets in the space of parameters $A+B$ and $C$:
\begin{itemize}
\item[ ] region I: $\{A+B-2C>0, A+B+2C>0\}$,
\item[ ] region II: $\{A+B-2C<0, A+B+2C>0\}$,
\item[ ] region III: $\{A+B-2C<0, A+B+2C<0\}$,
\item[ ] region IV: $\{A+B-2C>0, A+B+2C<0\}$.
\end{itemize}
Moreover, there are some non generic sections corresponding to the  
boundaries of the above regions, where at least one of the functions  
$A+B \pm 2C$ vanishes; in these case the corresponding paraboloid  
becomes a plane (an example is given by the section T). Figure~\ref{pic_Svari} describes the relation between the surfaces  
$S_1$, $S_2$ and $S_3$ for parameters belonging to region I.


\begin{Prop}
The rank of the distribution $\Delta_M$ is $6$ on $\Kmin-(S_1 \cup  
S_2)$.
\begin{Dim}
Since the determinant of the matrix (\ref{mat_min}) vanishes only on  
$S_1 \cup S_2$, outside of this set the rank of $\Delta_M$ is maximal.
\end{Dim}
\end{Prop}

\begin{Prop}
The rank of the distribution $\Delta_M$ on $S_1-S_2$ is~$5$.
\begin{Dim}
In order to study the rank of $\Delta_M$ on $S_1$, we set $\gamma=0$ in  
the matrix $M$ and calculate its adjoint matrix:
$$
\mathrm{adj}(M)|_{\gamma=0}=(\alpha^2-\beta^2 )\,
\left(
\begin{array}{cccccc}
0&0&0&0&0&0\\
0&0&0&0&0&0\\
0&0&0&0&0&0\\
0&0&0&0&0&0\\
0&0&0&0&0&0\\
\ast & \ast & 0 & 2(\beta^2-\alpha^2) & \ast & 2(\alpha^2-\beta^2)
\end{array}
\right)
$$
so the rank is lesser than 5 only when $\alpha^2=\beta^2$ that is on  
$S_1 \cap S_2$.
\end{Dim}
\end{Prop}

\begin{Prop}
The rank of the distribution $\Delta_M$ on $S_2-(S_1 \cup S_3)$ is~$5$.
\begin{Dim}
Let us study now the rank of the distribution on $S_2$ without its  
intersection with $S_1$: using the condition $\gamma\neq 0$ the  
equation (\ref{S2m}) of the two branches of $S_2$ becomes
\BAS
B = \frac{(\alpha-\beta)^2}{\gamma}-A+2C && \mbox{on the branch } B_1\\
B = \frac{(\alpha+\beta)^2}{\gamma}-A-2C && \mbox{on the branch } B_2\,.
\EAS
Substituting them in the matrix $M$ (\ref{mat_min}) and calculating  
the adjoint we obtain respectively
$$
\mathrm{adj}(M)|_{B_1}=(C\gamma-\alpha\beta)\left(
\begin{array}{cccccc}
\ast & \ast &-\gamma^2&\gamma^2&\ast &0 \\
\ast & \ast &-\gamma^2&\gamma^2&\ast &0 \\
\ast & \ast &2\gamma^2&-2\gamma^2&\ast &0 \\
\ast & \ast & \ast & \ast & \ast &0 \\
\ast & \ast & \ast & \ast & \ast &0 \\
\ast & \ast & \ast & \ast & \ast &0
\end {array}
\right)
$$
$$
\mathrm{adj}(M)|_{B_2}=(C\gamma-\alpha\beta)\left(
\begin{array}{cccccc}
\ast & \ast &-\gamma^2&-\gamma^2&\ast &0 \\
\ast & \ast &-\gamma^2&-\gamma^2&\ast &0 \\
\ast & \ast &-2\gamma^2&-2\gamma^2&\ast &0 \\
\ast & \ast & \ast & \ast & \ast &0 \\
\ast & \ast & \ast & \ast & \ast &0 \\
\ast & \ast & \ast & \ast & \ast &0
\end {array}
\right) \,.
$$
Then in both cases the rank is 5 except when $\gamma C=\alpha\beta$,  
that is outside of $S_3=B_1 \cap B_2$.
\end{Dim}
\end{Prop}
The surface $S_2-(S_1 \cup S_3)$ is formed by several connected  
components: in subsection~\ref{sezione_discreta} we will study which of  
these components are mapped one in the other by the discrete  
transformations and then generate the same type of coordinates system.

\begin{Prop}
The rank of the distribution $\Delta_M$ on $S_3-S_1$ is~$4$.
\begin{Dim}
 From the previous proposition the rank of $\Delta_M$ on $S_3-S_1$ is at  
most 4 and it is lesser if all the $4\times 4$ minors  of $M$ vanish.  
Outside of the intersection with $S_1$ (i.e. for $\gamma\neq0$) the  
equations (\ref{S3m}) of $S_3$ are
$$
B=\frac{\alpha^2+\beta^2}{\gamma}\,, \qquad C=\frac{\alpha\beta}{\gamma}
$$
and substituting these conditions in the matrix $M$ we obtains, by  
eliminating the second and third columns and the third and forth rows,  
the $4\times 4$ minor
$$
\left| \begin{array}{cccc}
0 &0 & -\gamma &0\\
-2\alpha & -\gamma & 0 &0\\
1 & 0 & 0 & 0\\
A & \alpha & \beta & \gamma
\end{array}\right| =\gamma^3\neq0\,.
$$
Hence, the rank of $\Delta_M$ on $S_3-S_1$ is always $4$ .
\end{Dim}
\end{Prop}

Let us now analyze the intersection between $S_1$ and $S_2$:  $S_1 \cap  
S_2$ is formed by two branches isomorphic to $\Reali^4$ intersecting in  
the three-dimensional vector space $\alpha=\beta=\gamma=0$. The first  
branch is described by the equations $\gamma=0$ and $\alpha=\beta$,  
while the second by the equations $\gamma=0$ and $\alpha=-\beta$.  
Inside $S_1 \cap S_2$ we point out the surface (union of two branches  
named, respectively, $C_1$ and $C_2$)
\BE\label{S4m}
S_4: \left\{
\begin{array}{l}
\gamma=0 \\
\alpha = \beta \\
A+B=2C
\end{array}
\right.
\;\bigcup \quad
\left\{
\begin{array}{l}
\gamma=0 \\
\alpha = -\beta \\
A+B=-2C
\end{array}
\right.
\qquad \mbox{dim }S_4=3
\EE
\begin{figure}
\centering
\includegraphics[width=6cm]{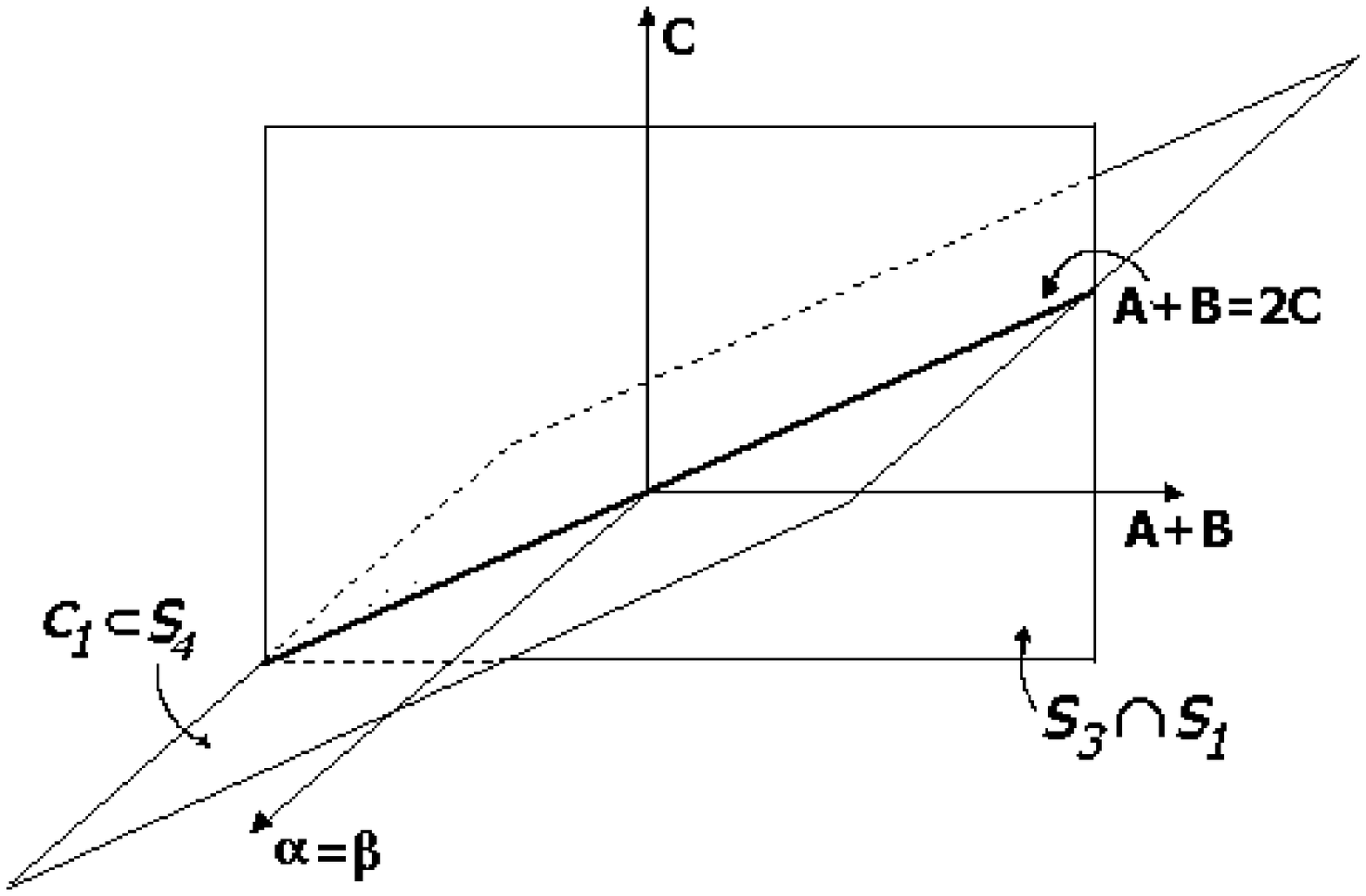}
\includegraphics[width=5cm]{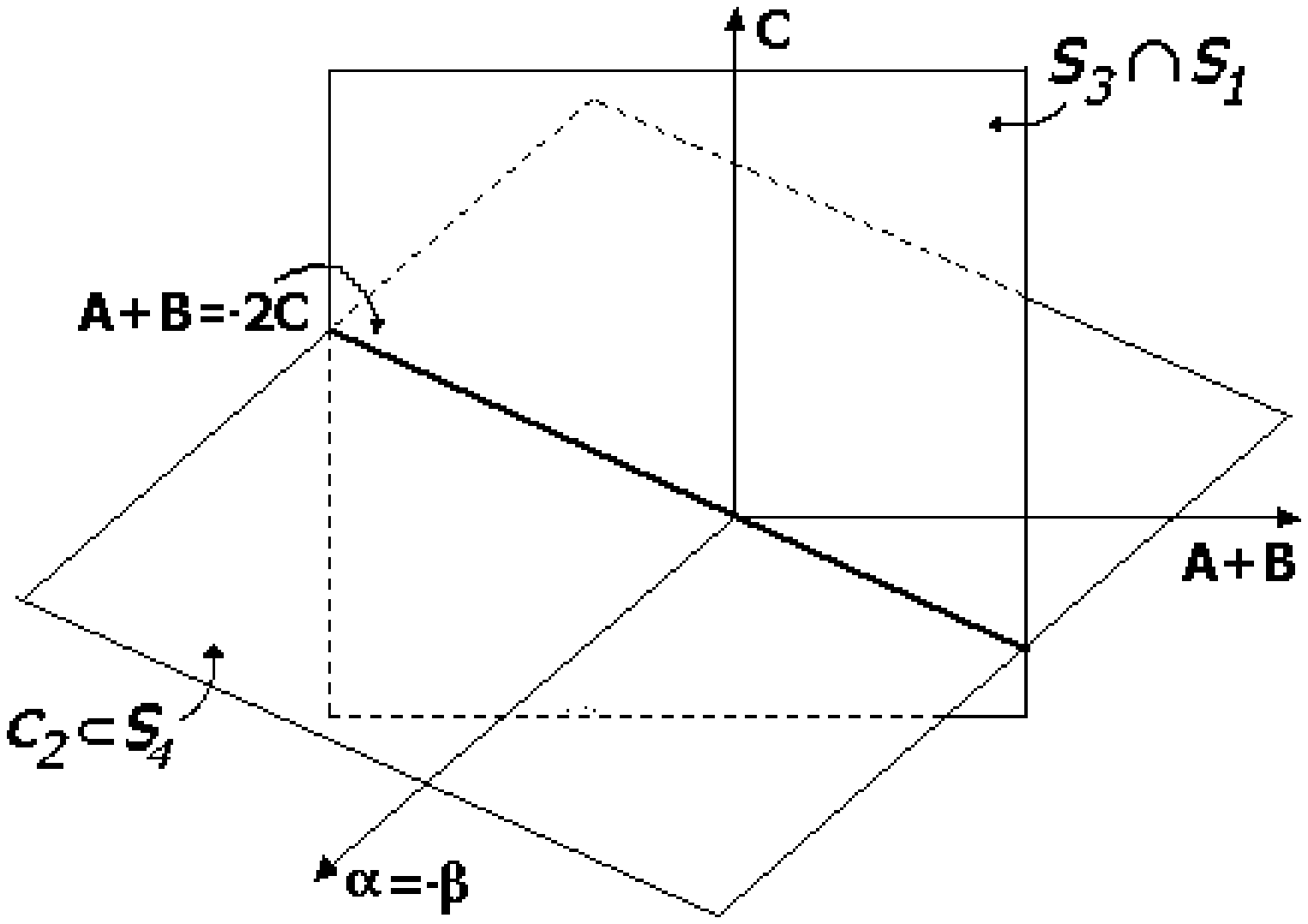}
\caption{The two branches of $S_4$ in $S_1\cap S_2$}\label{pic_S4}
\end{figure}
We observe that the branches $C_1$ and $C_2$ are both isomorphic to  
$\Reali^3$ and their intersection (belonging entirely to $S_1 \cap  
S_3$) is the line
\BE\label{S5m}
C_1 \cap C_2 = S_5: \left\{
\begin{array}{l}
\alpha=\beta=\gamma=0 \\
B = -A \\
C = 0
\end{array}
\right.
\qquad \mbox{dim }S_5=1
\EE
\begin{figure}
\centering
\includegraphics[width=6cm]{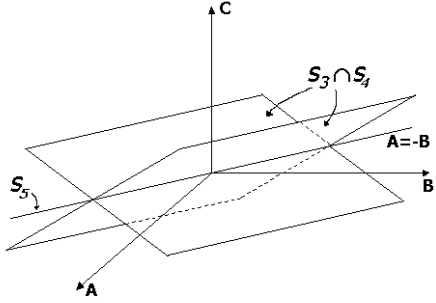}
\caption{The space $S_3\cap S_1$}\label{pic_S3}
\end{figure}

\begin{Prop}
The rank of the distribution $\Delta_M$ on $(S_1 \cap S_2)-(S_3 \cup  
S_4)$ is~$4$.
\begin{Dim}
We study the rank of $\Delta_M$ on the two branches of $S_1 \cap S_2$  
separately. Moreover we work outside of the intersection with $S_3$  
(that is we impose that both $\alpha$ and $\beta$ are different from  
zero). On the first branch, where $\gamma=0$ and $\alpha=\beta$, we get  
that all the non-vanishing $4\times 4$ minors are equal or proportional  
to $\alpha^2\,(2C-B-A)$; on the other branch, where $\gamma=0$ and  
$\alpha=-\beta$, all the non-vanishing minors are equal or proportional  
to $\alpha^2\,(B+A+2C)$. Then the rank is equal to $4$ outside of  
$S_4$.
\end{Dim}
\end{Prop}

\begin{Prop}
The rank of the distribution $\Delta_M$ on $S_4-S_3$ is~$3$.
\begin{Dim}
We study the rank of the distribution on the two branches $C_1$ and  
$C_2$ of $S_4$. On $C_1$, where $\beta=\alpha$ and $2C=A+B$ the vector  
fields $V_i$ form the matrix
$$
M|_{C_1}=\left(
\begin{array}{cccccc}
0&-2\,\alpha&-\alpha&0&0&0\\
-2\,\alpha&0&-\alpha&0&0&0\\
A+B&
A+B&A+B&\alpha&\alpha&0\\
2\,A&2\,B&A+B&\alpha&\alpha&0\\
1&-1&0&0&0&0\\
A&B&\frac{A+B}{2}&\alpha&\alpha&0
\end {array}
\right)
$$
and it is straightforward that
\BAS
V_2 &=& V_1-2\alpha V_5 \\
V_4 &=& V_3+(A-B)V_5 \\
2V_6 &=& V_3 +V_4 + \frac{A+B}{2\alpha}(V_1+V_2)
\EAS
and so outside of $S_4 \cap S_3$, where $\alpha=0$, the independent  
vector fields are $V_1$, $V_3$ and $V_5$ only. In a similar way, on  
$C_2$ $\beta=-\alpha$ and $2C=-A-B$ hold, and the vector fields $V_i$  
form the matrix
$$
M|_{C_2}=\left(
\begin{array}{cccccc}
0 & 2\alpha & -\alpha &  0 & 0 & 0 \\
-2\alpha & 0 &  \alpha &  0 & 0 & 0 \\
-A-B & -A-B &  A+B & -\alpha & \alpha & 0 \\
2A & 2B & -A-B & \alpha &  -\alpha & 0 \\
1 & -1 &  0 & 0 &  0 & 0 \\
A & B & -\frac{A+B}{2} & \alpha & -\alpha & 0
\end {array}
\right)
$$
and so
\BAS
V_2 &=& -V_1-2\alpha V_5 \\
V_4 &=& -V_3+(A-B)V_5 \\
2V_6 &=& -V_3 +V_4 + \frac{A+B}{2\alpha}(V_2-V_1)\,.
\EAS
The rank of the distribution on the two branches is then given by the  
rank of the two matrices
$$
\left(
\begin{array}{cccccc}
0&-2\,\alpha&-\alpha&0&0&0\\
A+B&A+B&A+B&\alpha&\alpha&0\\
1&-1&0&0&0&0
\end {array}
\right)
$$
$$
\left(
\begin{array}{cccccc}
0 & 2\alpha &  -\alpha &  0 & 0 &  0 \\
-A-B & -A-B &  A+B & -\alpha & \alpha & 0 \\
1 & -1 & 0 & 0 &  0 & 0
\end {array}
\right)\,.
$$
Because all the non-vanishing $3\times 3$ minors of these two matrices  
are proportional to $\alpha^2$, outside of  $S_4 \cap S_3$ the rank is  
3.
\end{Dim}
\end{Prop}

The rank of the distribution on $S_1 \cap S_3$ remains to be evaluated.  
The space $S_1 \cap S_3$ (see Figure~\ref{pic_S3}) is a  
three-dimensional vector space described by the equations  
$\alpha=\beta=\gamma=0$, with coordinates $A$, $B$ and $C$. We recall  
that \mbox{$S_3 \cap S_4 \subset S_3 \cap S_1$}. On $S_3 \cap S_1$ the  
only independent $V_i$ are $V_3$, $V_5$ and $V_6$ and their components,  
with respect to $(\partial_A,\partial_B,\partial_C)$ can be collected  
in the matrix
\BE\label{mat_min_rid}
\widetilde{M}=\left(
\begin{array}{ccc}
2C & 2C & A+B  \\
1 & -1 &  0 \\
A &  B &  C
\end{array}
\right)\,.
\EE

\begin{Prop}
The rank of the distribution $\Delta_M$ on $(S_3 \cap S_1)-S_4$ is~$3$.
\begin{Dim}
The determinant of the matrix (\ref{mat_min_rid}) is $\det  
\widetilde{M} =(A+B+2C)(2C-A-B)$ and it vanishes only on the  
intersection with $S_4$.
\end{Dim}
\end{Prop}

\begin{Prop}
The rank of the distribution $\Delta_M$ on $(S_3 \cap S_4)-S_5$ is~$2$.
\begin{Dim}
Evaluating the matrix $\widetilde{M}$ on the two branch of $S_3 \cap  
S_4$ one obtains two matrices whose adjoint have the form
$$
\mathrm{adj}(\widetilde{M})=(A+B)  \left(
\begin{array}{ccc}
1&-1&A-B\\
1&-1&A-B\\
-2&2&2(B-A)
\end {array}
\right)
$$
$$
\mathrm{adj}(\widetilde{M})=-(A+B)  \left(
\begin{array}{ccc}
  1&1&B-A\\
  1&1&B-A\\
  2&2&2(B-A)
\end {array}
  \right)
$$
and then the rank is lesser than 2 only on the intersection of the two  
branches given by $B=-A$, that is on the line $S_5$.
\end{Dim}
\end{Prop}

\begin{Prop}
The rank of the distribution $\Delta_M$ on $S_5$ is~$1$.
\begin{Dim}
On $S_5$ the only independent vector field is the constant vector $V_5$.
\end{Dim}
\end{Prop}

\subsection{Discrete transformations}\label{sezione_discreta}
As we already mentioned besides the continuous transformations  
associated with the vector fields $V_i$ we have to consider also some  
discrete transformation leaving unchanged the web associated with a  
given Killing tensor: the first one is the change of the sign of the  
tensor
$$
R_0:K \to -K
$$
and the others are induced from the discrete isometries of the  
Minkowski plane $\{ \bar t=t,\;\bar x=-x\}$
and $\{\bar t=-t,\;\bar x=x\}$, they are
$$
R_1: C\to-C, \alpha\to-\alpha
$$
$$
R_2: C\to-C, \beta\to-\beta\,.
$$
Now we have to study the connected components of the sets determined in  
the subsection~\ref{sezione_continua} and which of them are linked  
through one of the above discrete transformations. Since some of these  
sets have a quite high dimension we use the sectioning technique  
presented in subsection~\ref{sezione_inco}.

In order to study the open set $\Kmin-S_1-S_2$ we observe that it is  
the set where the three functions
\BAS
\gamma \\
Z_+ &=& \gamma(A+B-2C)-(\alpha-\beta)^2 \\
Z_- &=& \gamma(A+B+2C)-(\alpha+\beta)^2\\
\EAS
are all different from zero, where the notation of~\cite{McL} has been used.
 Then, the continuous function  
$\Kmin\to\Reali^3$ given by $\Phi=(\gamma,Z_+,Z_-)$ maps  
$\Kmin-S_1-S_2$ in the eight connected components of $\Reali^3$ without  
the coordinate planes. We introduce the eight (not empty) sets  
$\Gamma_1,\ldots,\Gamma_8$ such that
$$
\begin{array}{cc}
\Phi(\Gamma_1)=(+,+,+), & \Phi(\Gamma_5)=(-,+,+), \\
\Phi(\Gamma_2)=(+,+,-), & \Phi(\Gamma_6)=(-,+,-), \\
\Phi(\Gamma_3)=(+,-,+), & \Phi(\Gamma_7)=(-,-,+), \\
\Phi(\Gamma_4)=(+,-,-), & \Phi(\Gamma_8)=(-,-,-).
\end{array}
$$
The sets $\Gamma_1,\ldots,\Gamma_8$ form a partition of $\Kmin-S_1-S_2$.

\begin{Prop}
All the sets $\Gamma_1, \ldots,\Gamma_8$ are connected and then the set  
$\Kmin-(S_1 \cup S_2)$ has $8$ connected components. We have three  
orbits of the action: $\Gamma_1 \cup \Gamma_5$, $\Gamma_4 \cup  
\Gamma_8$  (both linked by the transformation $R_0$) and $\Gamma_2 \cup  
\Gamma_3 \cup \Gamma_6 \cup \Gamma_7$ (linked by $R_0$ and $R_1$).
\begin{Dim}
We prove that all the sets $\Gamma_1, \ldots,\Gamma_8$ are connected by  
showing that all the connected components of their sections are  
equivalent in the sense of Lemma~\ref{inco}. First of all we remark  
that for any two sections with parameter belonging to same region of  
Figure~\ref{pic_S2}, the corresponding connected components (of any  
$\Gamma_i$) are trivially equivalent. Hence the sets $\Gamma_1$ and  
$\Gamma_5$ are connected because their sections are not empty only for  
parameters in regions I and III, respectively. For any of the other  
$\Gamma_i$ there exists a region of the parameter space in which the  
corresponding section has a unique connected component. On the other  
hand it is possible to construct a continuous path connecting a point  
of any connected components of the section of a given $\Gamma_i$ to a  
point in the one with a unique connected component. For instance moving  
the parameters along the path shown in Figure~\ref{pic_S2} (with  
$A+B+2C$ constant) and leaving $\alpha$, $\beta$ and $\gamma$ fixed we  
link any point in one of the two connected components of $\Gamma_3$  
obtained for parameter in region I, to a point in the unique  
connected component of the section of $\Gamma_3$, obtained for  
parameters in region II.  We conclude, by applying Lemma~\ref{inco},  
that the sets $\Gamma_i$ are all connected. From the definition of  
transformation $R_0$ and $R_1$ we get that $\Gamma_1$ is in the same  
orbit of $\Gamma_5$, $\Gamma_4$ is in the same orbit of $\Gamma_8$ and  
all the other $\Gamma_i$ are in the same orbit. Moreover, because all  
the transformations $R_i$ maps any of the sets $\Gamma_1 \cup  
\Gamma_5$, $\Gamma_4 \cup \Gamma_8$, $\Gamma_2 \cup \Gamma_3 \cup  
\Gamma_6 \cup \Gamma_7$ into themselves, they form distinct orbits.
\end{Dim}
\end{Prop}

In order to study the set $S_2 - (S_1 \cup S_3)$ we introduce the eight  
(not empty) sets $\Theta_1,\ldots,\Theta_8$ such that
$$
\begin{array}{cc}
\Phi(\Theta_1)=(+,0,+), & \Phi(\Theta_5)=(-,0,+), \\
\Phi(\Theta_2)=(+,0,-), & \Phi(\Theta_6)=(-,0,-), \\
\Phi(\Theta_3)=(+,+,0), & \Phi(\Theta_7)=(-,+,0), \\
\Phi(\Theta_4)=(+,-,0), & \Phi(\Theta_8)=(-,-,0).
\end{array}
$$
The sets $\Theta_1,\ldots,\Theta_8$ form a partition of $S_2 - (S_1  
\cup S_3)$.

\begin{Prop}
All the sets $\Theta_1,\ldots,\Theta_8$ are connected and then the set  
$S_2 - (S_1 \cup S_3)$ has $8$ connected components. We have two orbits  
of the action: $\Theta_1 \cup \Theta_3 \cup \Theta_5 \cup \Theta_7$ and  
$\Theta_2 \cup \Theta_4 \cup \Theta_6 \cup \Theta_8$ (both linked by  
the transformations $R_0$ and $R_1$).
\begin{Dim}
As in the previous proposition we prove that all the sets $\Theta_1,  
\ldots,\Theta_8$ are connected by showing that all the connected  
components of their sections are equivalent in the sense of  
Lemma~\ref{inco}. Also in this case for any two sections with parameter  
belonging to same region of Figure~\ref{pic_S2}, the corresponding  
connected components (of any $\Theta_i$) are trivially equivalent. We  
observe that because the plane $\gamma=0$ is removed all the  
paraboloids that form the sections of the sets $\Theta_i$ consist of  
at least two disconnected parts. For any $\Theta_i$ a section with a  
unique connected component exists. For instance in the section labeled
T, in Figure~\ref{pic_S2}, the sets $\Theta_1$ and $\Theta_5$ have a  
connected section. Moreover it is possible to construct a continuous  
path connecting a point of any connected components of the section of a  
given $\Theta_i$ to a point in the one with a unique connected  
component. We conclude, by applying Lemma~\ref{inco}, that the sets  
$\Theta_i$ are all connected. From the definition of transformation  
$R_0$ and $R_1$ we get that $\Theta_1$,  $\Theta_3$, $\Theta_5$ and  
$\Theta_7$ are mapped one into the other and then are in the same  
orbit, as well $\Theta_2$,  $\Theta_4$, $\Theta_6$ and $\Theta_8$.  
Moreover, because all the transformations $R_i$ maps the two sets  
$\Theta_1 \cup \Theta_3 \cup \Theta_5 \cup \Theta_7$ and $\Theta_2 \cup  
\Theta_4 \cup \Theta_6 \cup \Theta_8$ into themselves, they form  
distinct orbits.
\end{Dim}
\end{Prop}

\begin{Prop}
The set $S_3 - S_1$ has two connected components mapped one in the  
other by $R_0$, and so it forms a unique orbit of the action.
\begin{Dim}
A section of the set $S_3 - S_1$ is not empty only if its parameters  
(A,B,C) belong to the closure of regions I and III referring to the  
notation of Figure~\ref{pic_S2}). All sections with parameters  
belonging to the interior of I (respectively, III) have four connected  
components which are equivalent to the positive part of the axis  
$\gamma$ (respectively, the negative part) in the sections with  
parameters $A+B=C=0$. The sections with parameters on the boundary of I  
(respectively, III) and $C\neq 0$ have two connected components, also  
equivalent to the positive part of the axis $\gamma$ (respectively, the  
negative part) in the sections with parameters $A+B=C=0$. Hence just  
two different equivalence classes of sections exist, corresponding to  
two connected components of $S_3 - S_1$: one corresponds to positive  
values of $\gamma$, while the other to negative ones. Since the  
transformation $R_0$ maps the positive part of the $\gamma$ axis in the  
negative one, these two connected components form a unique orbit.
\end{Dim}
\end{Prop}

\begin{Prop}
The set $S_1 - S_2$ is is formed by $4$ connected components. Each pair of components symmetric with respect to the origin are linked by $R_0$, thus two orbit of the action are present.
\begin{Dim}
$S_1$ is homeomorphic to $\Reali^5$, then it is divided in four  
connected parts by the two four-dimensional hyperplanes that form $S_1  
\cap S_2$. The transformation $R_0$ represents a central symmetry and links together components symmetric with respect to the origin. The two transformations $R_1$ and $R_2$ on $S_1-S_2$ are symmetries with respect to the $\alpha$ and $\beta$ axes, hence they map the two orbits into themselves. 
\end{Dim}
\end{Prop}

\begin{Prop}
The set $(S_1 \cap S_2) - (S_3 \cup S_4)$ contains $8$ connected  
components. They can be linked together using the three discrete  
transformation $R_0$, $R_1$ and $R_2$ and so they form a unique orbit  
of the action.
\begin{Dim}
Each of the two branches of $S_1 \cap S_2$ is homeomorphic to  
$\Reali^4$. Cutting out from the first one the two three-dimensional  
hyperplanes $S_3 \cap S_1$ and $C_1$, and  from the second one $S_3  
\cap S_1$ and $C_2$, respectively, we obtain four connected components  
in each case. In order to prove that all the components form a unique  
orbit, we can consider just the branch defined by $\alpha=\beta$, because $R_1$ (or equivalently $R_2$) maps  
each branch into the other. The  
transformation $R_0$ maps one in the other the components symmetric  
with respect to the origin (see Figure~\ref{pic_S4}), while the  
transformation $R_1\circ R_2$ (corresponding in this branch to the inversion of the $\alpha=\beta$ axis) maps one in the other the components symmetric with  
respect to the hyperplane $S_3 \cap S_1$. Therefore all the eight connected  
components form a unique orbit.
\end{Dim}
\end{Prop}

\begin{Prop}
The set $S_4 - S_3$ contains $4$ connected components. They can be  
linked together using $R_0$ and one between $R_1$ and $R_2$ and so they  
form a unique orbit of the action.
\begin{Dim}
Each of the two branches $C_1$ and $C_2$ of $S_4$ is homeomorphic to  
$\Reali^3$. Cutting out from the first one the two two-dimensional  
hyperplane $S_3 \cap S_4$, we obtain two connected components in each  
case. Each of the transformations $R_1$ or $R_2$ maps $C_1$ in $C_2$ and $R_0$ links the two connected components of $C_2$. Thus we have a unique  
orbit.
\end{Dim}
\end{Prop}

\begin{Prop}
The space $(S_3 \cap S_1) - S_4$ has four connected components. The two  
pairs of components symmetric with respect to the origin (linked by  
$R_0$) form two different orbits of the action.
\begin{Dim}
As shown in Figure~\ref{pic_S3}, in the space of coordinates $A$, $B$  
and $C$, the set $(S_3 \cap S_1) - S_4$ is composed by the four dihedra  
determined by the two planes $A+B=2C$ and $A+B=-2C$. Hence it has four  
connected components. The transformation $R_0$ links together the two  
dihedra containing the plane $C=0$ as well as the other pair of  
dihedra. Unfortunately neither $R_1$ nor $R_2$ is able to connect  
together these two pairs of dihedra. Hence in $(S_3 \cap S_1) - S_4$ we  
have two different orbits: indeed, the KT's belonging to the  
pair that contains the plane $C=0$ define pseudo-Cartesian  
coordinates, while the ones belonging to the other pair are not  
characteristic tensors, with everywhere imaginary eigenvalues.
\end{Dim}
\end{Prop}

\begin{Prop}
The space $(S_3 \cap S_4) - S_5$ has four connected components. They  
are mapped one into the other by the two discrete transformations $R_0$  
and $R_1$, hence they form a unique orbit of the action.
\begin{Dim}
As shown in Figure~\ref{pic_S3}, in the space of coordinates $A$, $B$  
and $C$, the set $(S_3 \cap S_4) - S_5$ is formed by two planes  
intersecting on the line $A+B=0$, $C=0$ ($S_5$) without their  
intersection. Hence it has four connected components. The  
transformation $R_0$ maps an half of each plane in the other; moreover,  
the transformation $R_1$  maps each plane into the other.
\end{Dim}
\end{Prop}

\begin{Prop}
The line $S_5$ formed by the (non-characteristic) tensors of the kind  
$\tau\,\Tens{g}$ is a connected orbit of the action.
\end{Prop}

The following list contains all orbits of the action of the group generated by the vector  
fields $V_i$, extended with the three finite transformations. For each orbit a representative tensor is given. Orbits of characteristic Killing tensors are labeled according both to \cite{McL} and \cite{Kalnins} and the associated \emph{complete} web is plotted. In each picture set of singular points and and the two distinct foliations of the web are amphasize completing the partial representation given in \cite{McL} and \cite{Kalnins}: the leaves belonging to the two foliations are plotted, respectively, dashed and continuous and the grey lines represent the boundaries of the singular set of the web).
\begin{description}
\item{M1)} The set $\Gamma_1 \cup \Gamma_5$, contained in $\Kmin - (S_1 \cup S_2)$,  
where $Z_+$ and $Z_-$ are both positive: SC9, elliptic coordinates of type I. A tensor of this type is:
$$
\left(
\begin{array}{cc}
x^2 & xt  \\
xt & t^2+1
\end{array}
\right)\,.
$$
\item{M2)} The set $\Gamma_2 \cup \Gamma_3 \cup \Gamma_6 \cup \Gamma_7$, contained  
in $\Kmin - (S_1 \cup S_2)$, where $Z_+$ and $Z_-$ have different sign:  
SC8, hyperbolic coordinates of type I. A tensor of this type is:
$$
\left(
\begin{array}{cc}
x^2 & 1+xt  \\
1+xt  & t^2
\end{array}
\right)\,.
$$
\item{M3)} The set $\Gamma_4 \cup \Gamma_8$, contained in $\Kmin - (S_1 \cup S_2)$,  
where $Z_+$ and $Z_-$ are both negative: SC5 and SC10, elliptic coordinates of type II. A tensor of this type is:
$$
\left(
\begin{array}{cc}
x^2 & xt  \\
xt & t^2-1
\end{array}
\right)\,.
$$
\item{M4)} The set $\Theta_1 \cup \Theta_3 \cup \Theta_5 \cup \Theta_7$ contained in  
$S_2-(S_1 \cup S_3)$, where the non-vanishing one of the two functions  
$Z_\pm$ is positive: SC6, hyperbolic coordinates of type II. Two tensors of this type are:
$$
\left(
\begin{array}{cc}
x^2+1 & xt+1  \\
xt+1 & t^2+1
\end{array}
\right)\,,\quad
\left(
\begin{array}{cc}
x^2+1 & xt-1  \\
xt-1 & t^2+1
\end{array}
\right)\,.
$$
\item{M5)} The set $\Theta_2 \cup \Theta_4 \cup \Theta_6 \cup \Theta_8$ contained in  
$S_2-(S_1 \cup S_3)$, where the non-vanishing one of the two functions  
$Z_\pm$ is negative: SC7, hyperbolic coordinates of type III. Two tensors of this type are:
$$
\left(
\begin{array}{cc}
x^2-1 & xt-1  \\
xt-1 & t^2-1
\end{array}
\right)\,,\quad
\left(
\begin{array}{cc}
x^2-1 & xt+1  \\
xt+1 & t^2-1
\end{array}
\right)\,.
$$
\item{M6)} The set $S_3-S_1$: SC2, polar coordinates. A tensor of this type is:
$$
\left(
\begin{array}{cc}
x^2 & xt  \\
xt & t^2
\end{array}
\right)\,.
$$
\item{M7)} The subset of $S_1 -S_2$ containing the $\alpha$ axis: first web for SC4, parabolic coordinate of type I. A tensor of this type is:
$$
\left(
\begin{array}{cc}
2x & t  \\
t & 0
\end{array}
\right)\,.
$$
\item{M8)} The subset of $S_1 -S_2$ containing the $\beta$ axis: second web for SC4, parabolic coordinate of type I. A tensor of this type is:
$$
\left(
\begin{array}{cc}
0 & x  \\
x & 2t
\end{array}
\right)\,.
$$
\item{M9)} The set $(S_1 \cap S_2)- (S_3 \cup S_4)$: SC3, parabolic coordinate of type II. Two tensors of this type are:
$$
\left(
\begin{array}{cc}
1+2x & x+t  \\
x+t & 1+2t
\end{array}
\right)\,,\quad
\left(
\begin{array}{cc}
2x+1 & x+t-1  \\
x+t-1 & 2t+1
\end{array}
\right)\,.
$$
\item{M10)} The set $S_4-S_3$: no characteristic tensors. A tensor of this type is:
$$
\left(
\begin{array}{cc}
2x & x+t  \\
x+t & 2t
\end{array}
\right)\,.
$$
\item{M11)} The subset of $(S_1 \cap S_3)-S_4$ containing the plane $C=0$: SC1, Cartesian coordinates. A tensor of this type is:
$$
\left(
\begin{array}{cc}
1 & 0  \\
0 & 0
\end{array}
\right)\,.
$$
\item{M12)} The subset of $(S_1 \cap S_3)-S_4$  not containing the plane $C=0$:  
no characteristic tensors. A tensor of this type is:
$$
\left(
\begin{array}{cc}
0 & 1  \\
1 & 0
\end{array}
\right)\,.
$$
\item{M13)} The set $(S_3 \cap S_4)-S_5$: no characteristic tensors. A tensor of this type~is:
$$
\left(
\begin{array}{cc}
1 & 1  \\
1 & 1
\end{array}
\right)\,.
$$
\item{M14)} The line $S_5$, containing tensors multiple of the metric.
\end{description}

\begin{center}
\begin{picture}(160,170)
\put(-15,0){\includegraphics[width=6cm]{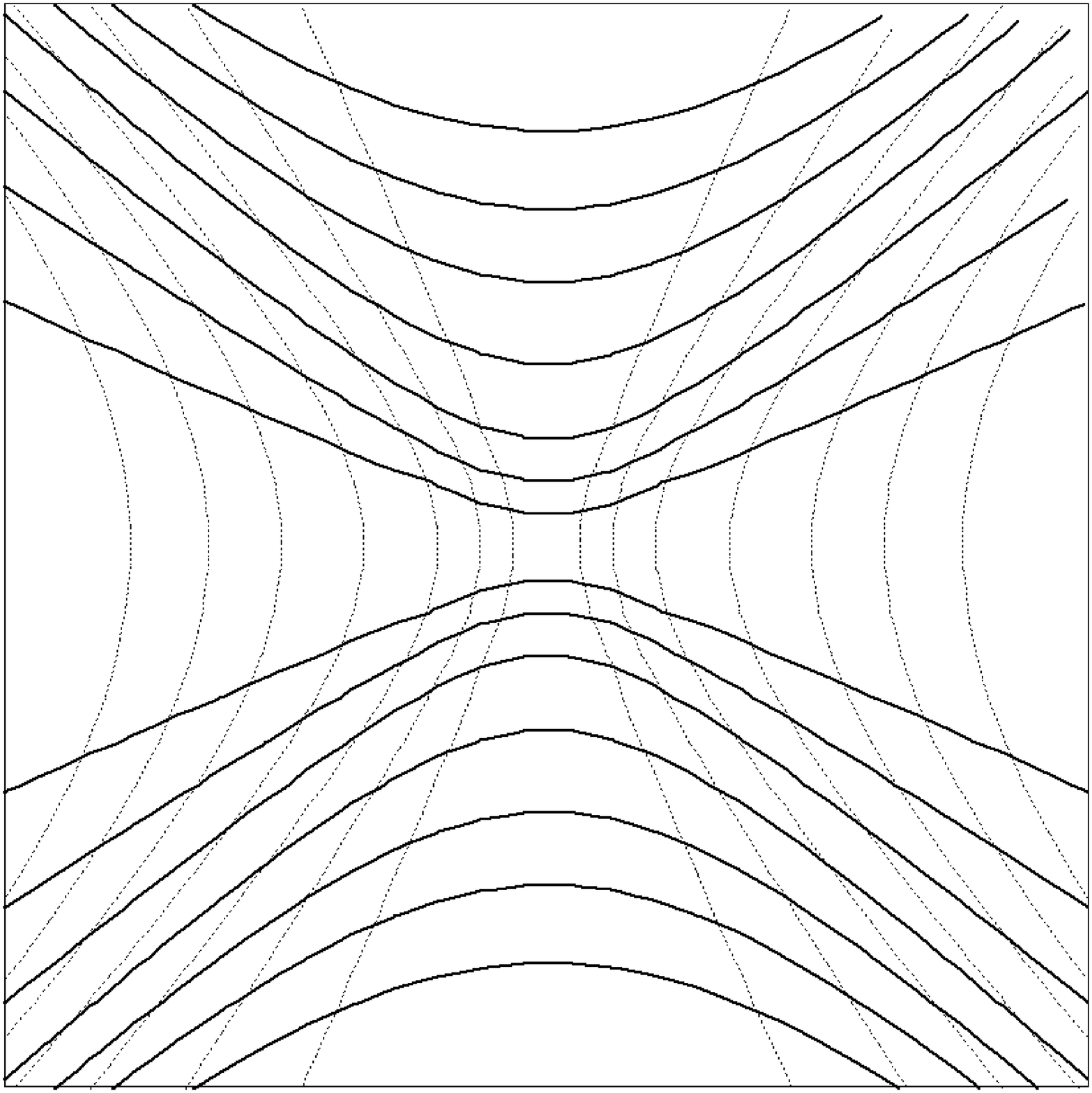}}
\put(40,0){Web for SC9}
\end{picture}
\begin{picture}(160,170)
\put(-15,0){\includegraphics[width=6cm]{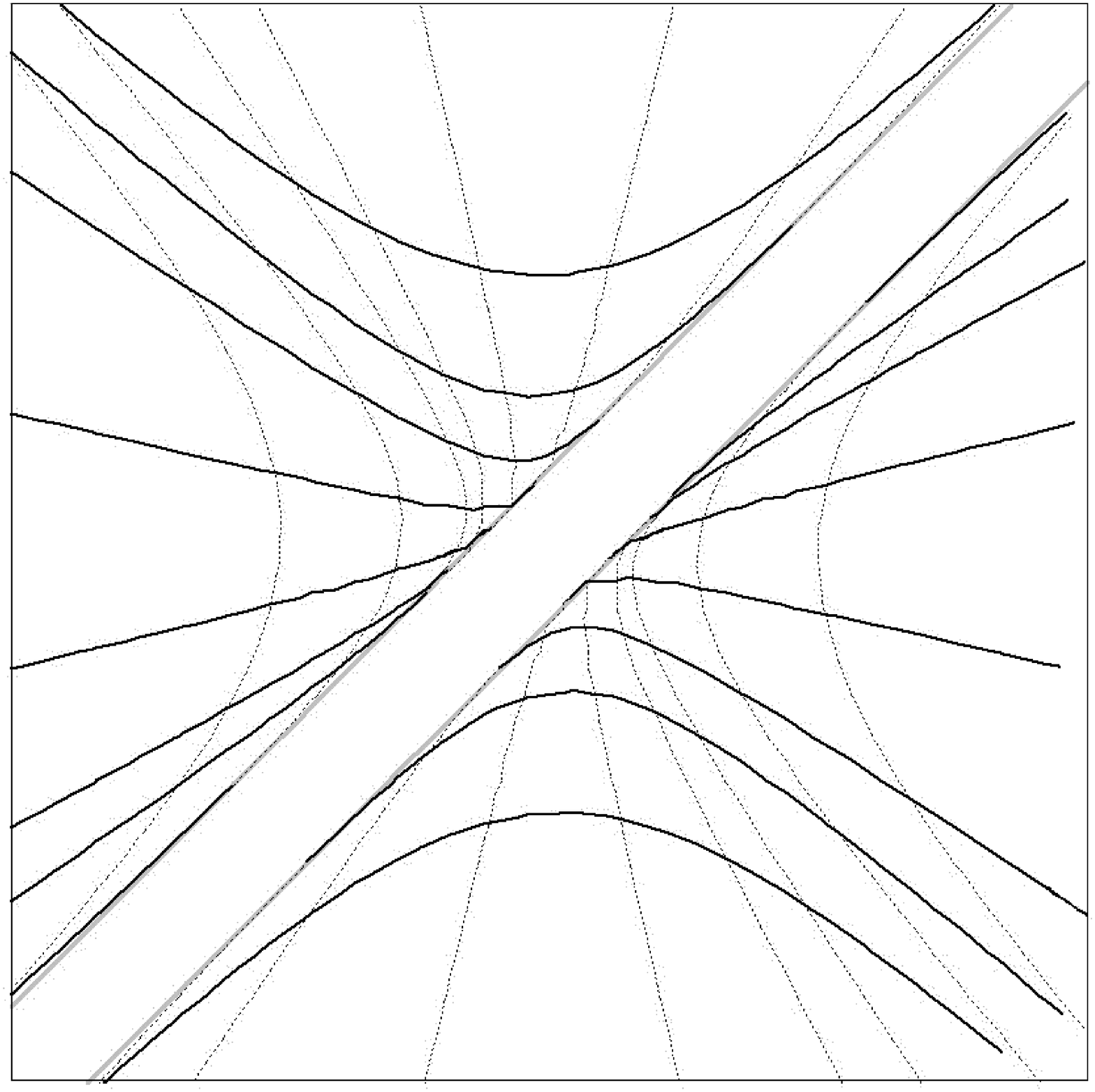}}
\put(40,0){Web for SC8}
\end{picture}
\end{center}
\begin{center}
\begin{picture}(160,170)
\put(-15,0){\includegraphics[width=6cm]{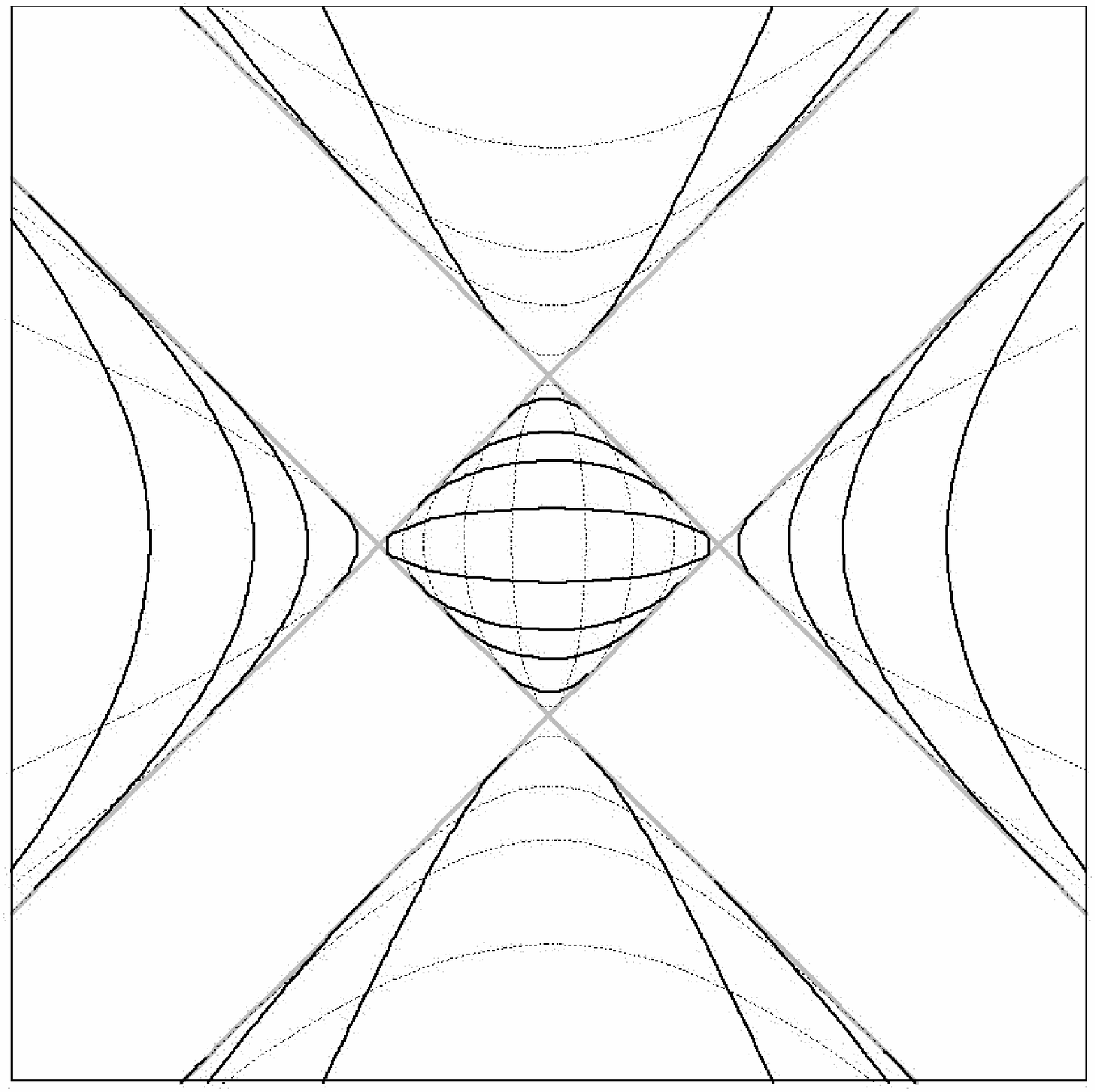}}
\put(20,0){Web for SC5 and SC10}
\end{picture}
\begin{picture}(160,170)
\put(-15,0){\includegraphics[width=6cm]{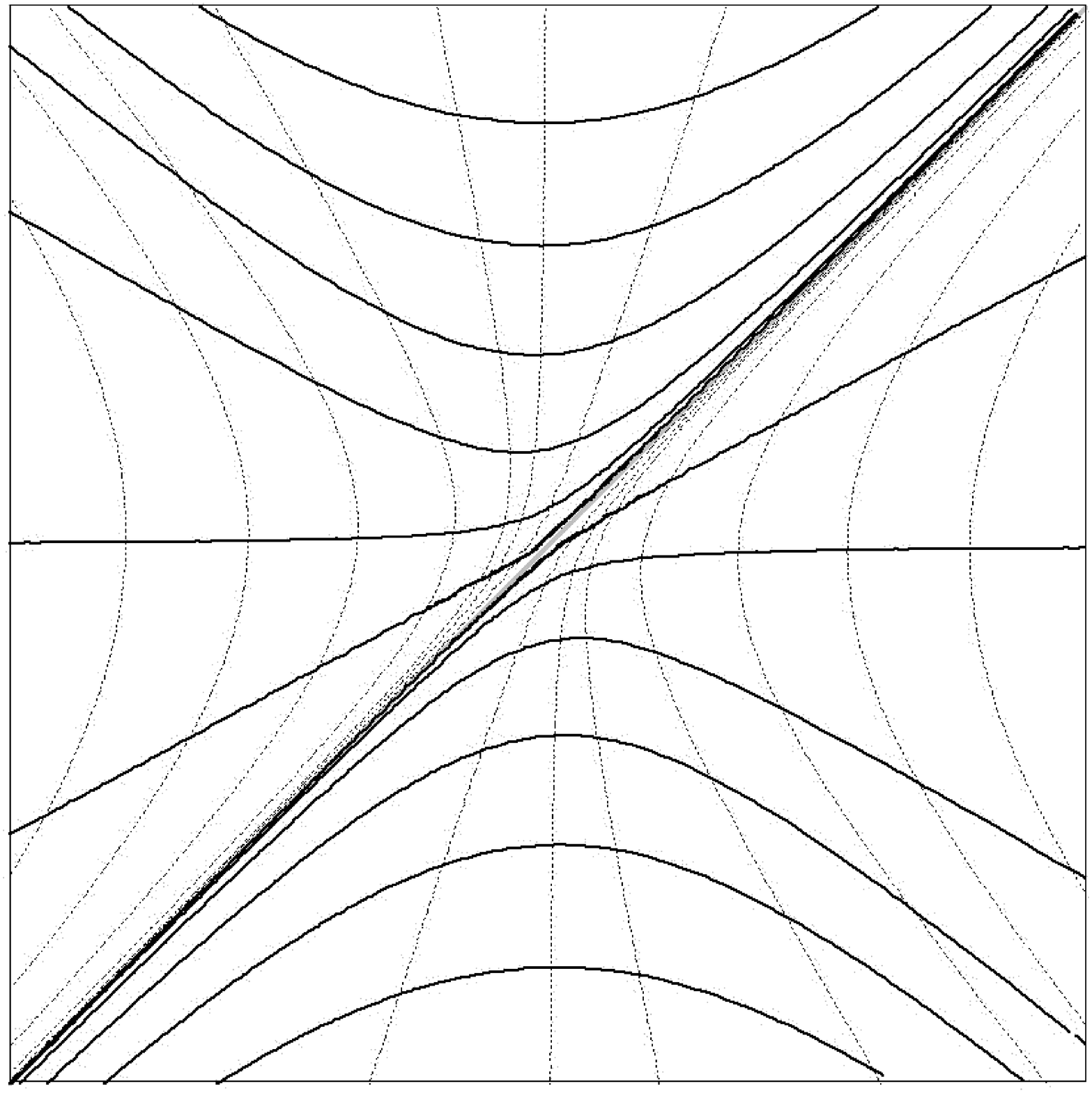}}
\put(40,0){Web for SC6}
\end{picture}
\end{center}
\begin{center}
\begin{picture}(160,170)
\put(-15,0){\includegraphics[width=6cm]{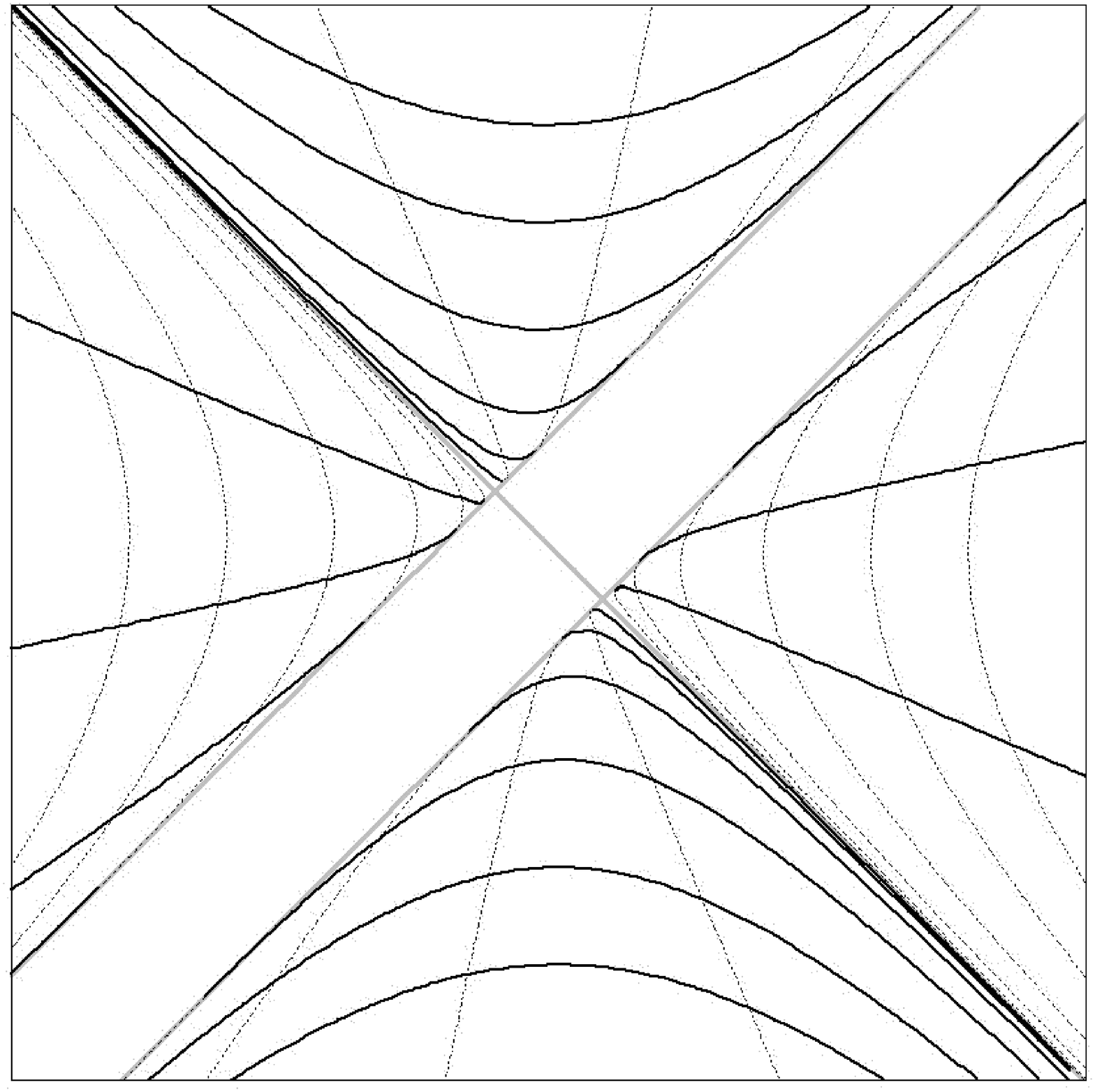}}
\put(40,0){Web for SC7}
\end{picture}
\begin{picture}(160,170)
\put(-15,0){\includegraphics[width=6cm]{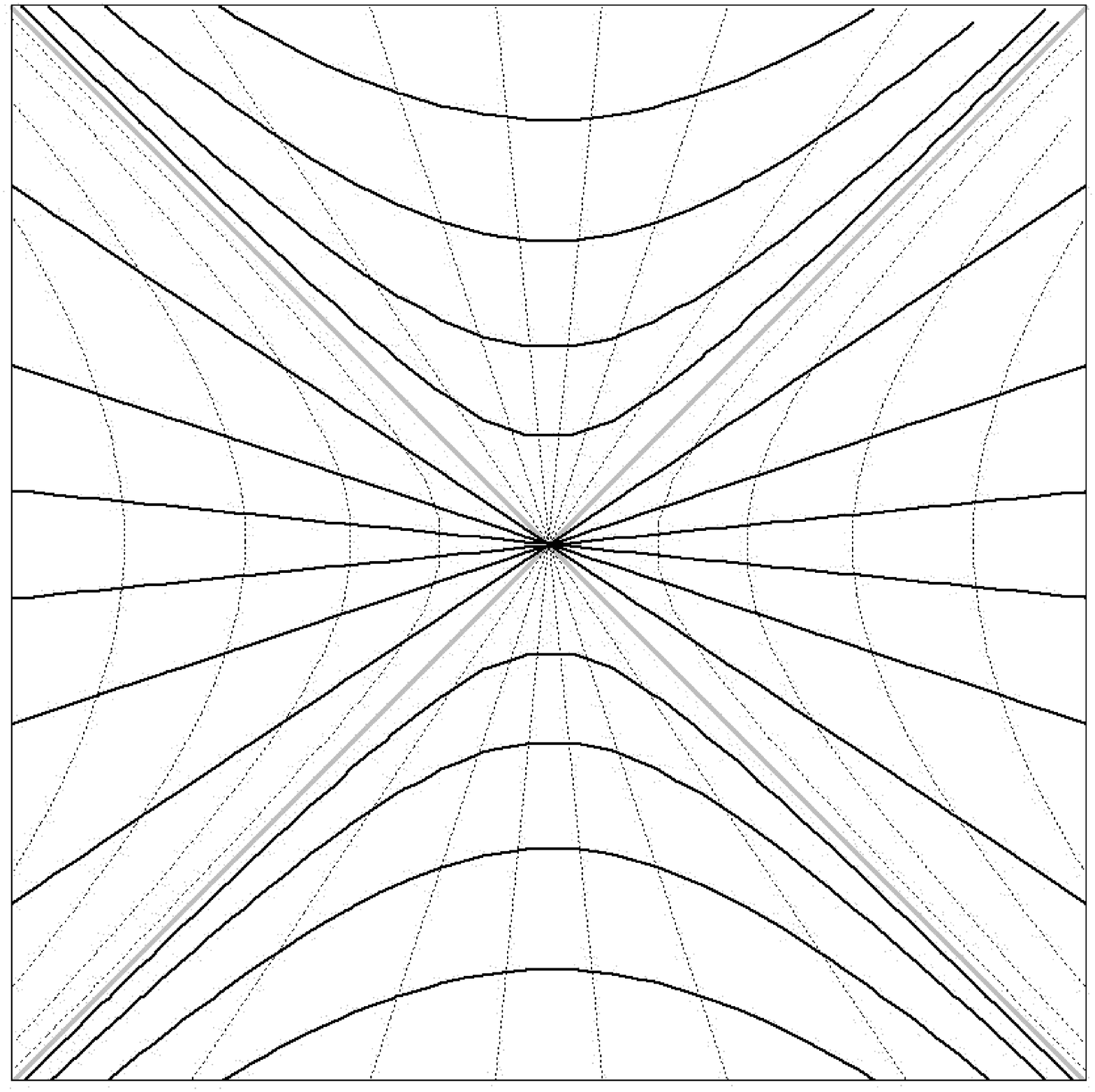}}
\put(40,0){Web for SC2}
\end{picture}
\end{center}
\begin{center}
\begin{picture}(160,170)
\put(-15,0){\includegraphics[width=6cm]{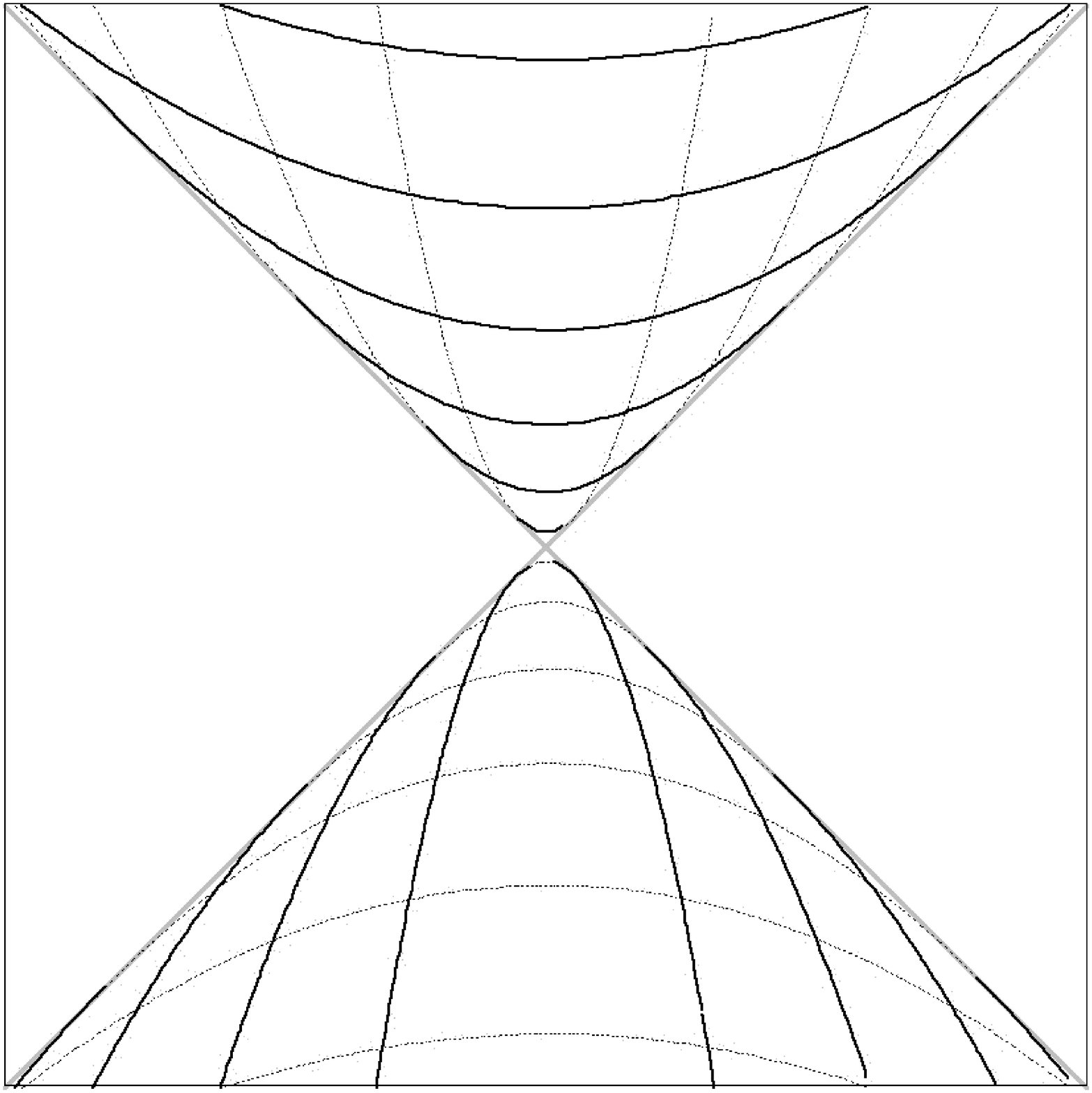}}
\put(30,0){First web for SC4}
\end{picture}
\begin{picture}(160,160)
\put(-15,0){\includegraphics[width=6cm]{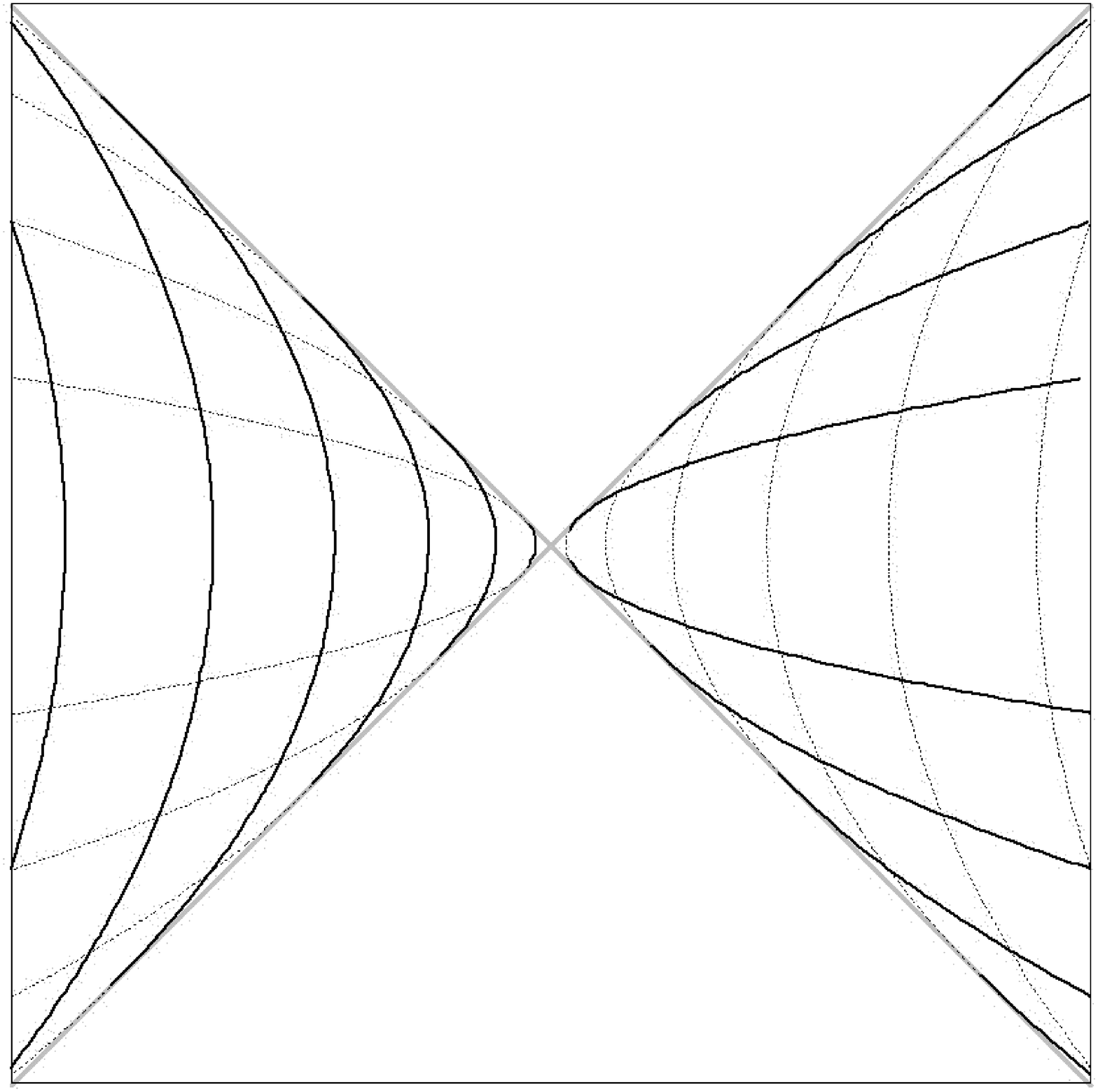}}
\put(25,0){Second web for SC4}
\end{picture}
\end{center}
\begin{center}
\begin{picture}(160,170)
\put(-15,0){\includegraphics[width=6cm]{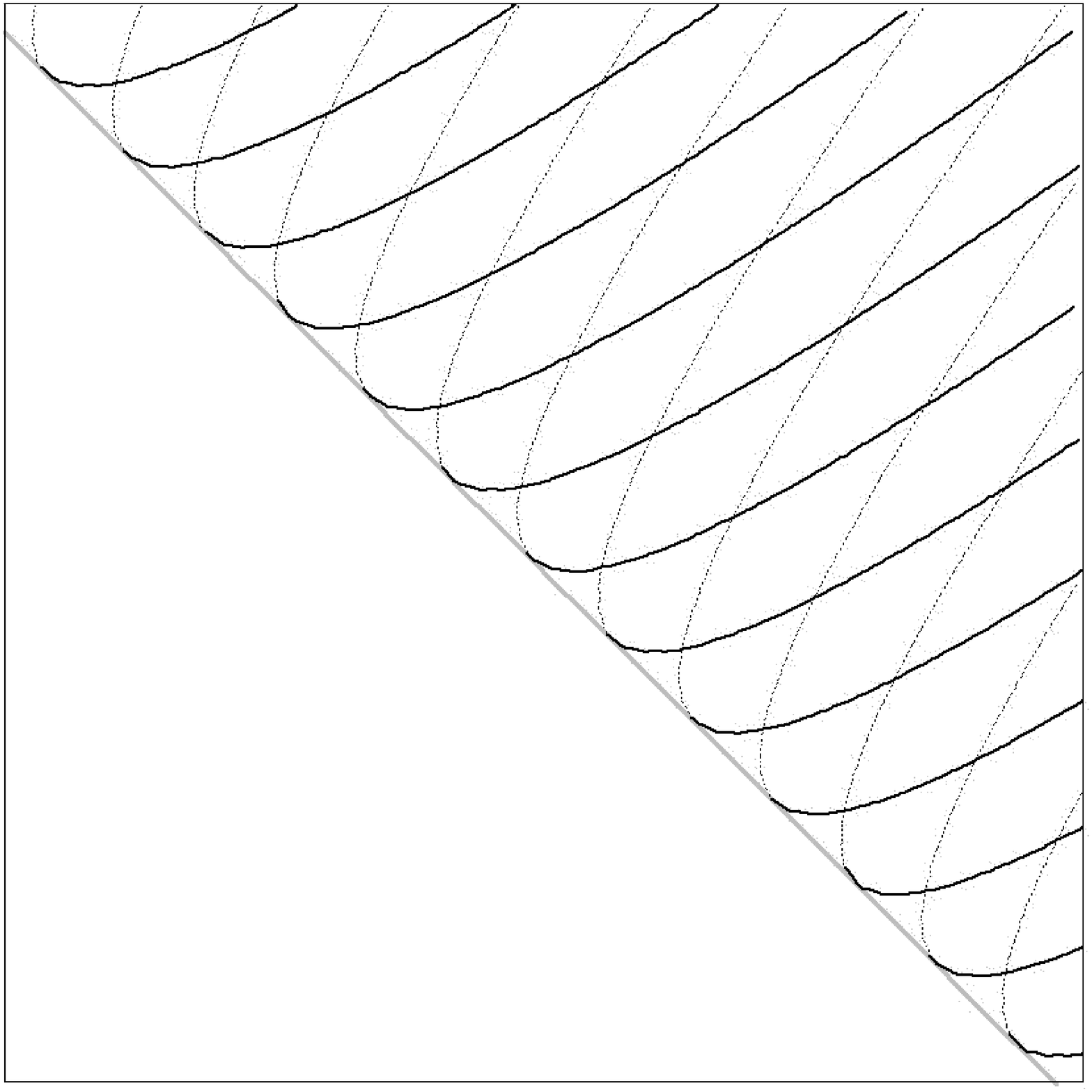}}
\put(40,0){Web for SC3}
\end{picture}
\begin{picture}(160,170)
\put(-15,0){\includegraphics[width=6cm]{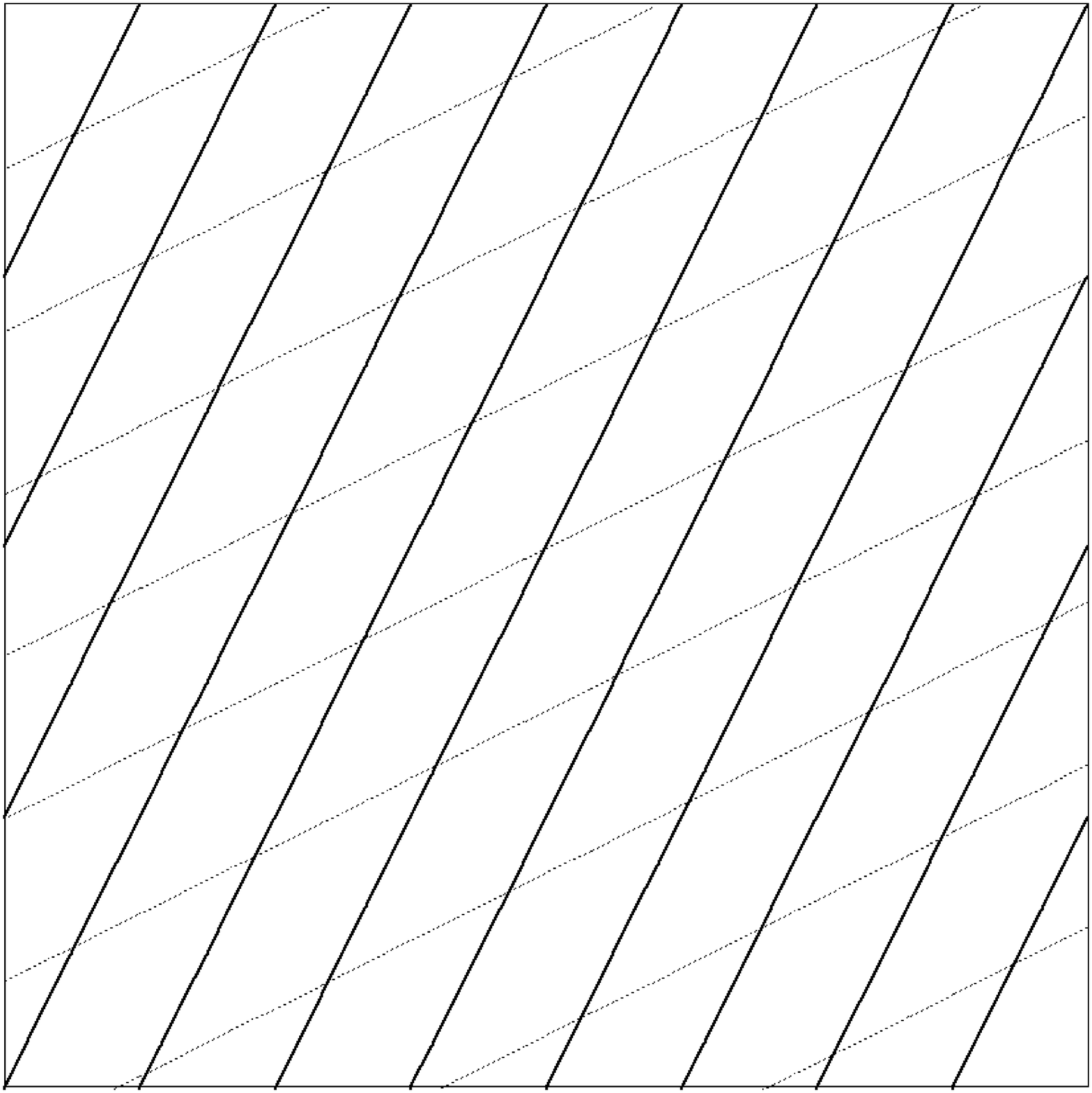}}
\put(40,0){Web for SC1}
\end{picture}
\end{center}

As in the Euclidean case, our classification is closely related to the one of Rastelli \cite{ref_9} based on 
the analysis of the singular set of the tensors.
The discriminant of the characteristic polynomial of the general KT of the Minkowski plane is
$$
\Delta= (\gamma(x+t)^2+2(\alpha+\beta)(x+t)+A+B+2C)(\gamma(x-t)^2+2(\alpha-\beta)(x-t)+A+B-2C).
$$ 
For $\gamma \neq 0$ (i.e.\ outside of $S_1$), we rewrite $\Delta$ as
$$
\gamma^2\left( (x+t+\textstyle\frac{\alpha+\beta}{\gamma})^2+\frac{1}{\gamma^2}Z_+
\right)
\left( \big(x-t+\textstyle\frac{\alpha-\beta}{\gamma}\big)^2+\frac{1}{\gamma^2}Z_- \right)\,.
$$
In this case the set $\Delta=0$ is made of two couples of lines parallel to $x=t$ and $x=-t$, respectively.
It is immediate to see that the lines of the first (second) pair are real and distinct, real and coinciding, imaginary
according to the fact that $Z_+$ ($Z_-$) is negative, zero or positive.
So the singular set is empty when both $Z_\pm$ are positive (SC9); a strip when $Z_+Z_- <0$ (SC8); two intersecting strips without their intersection when $Z_\pm$ are negative (SC5, SC10); a line when one of $Z_\pm$ vanishes and the other is positive (SC6); a strip and a line orthogonal to it when one of $Z_\pm$ vanishes and the other is positive (SC7); two orthogonal lines if both $Z_\pm$ vanish. 

On $S_1$ we have $\gamma=0$ and the discriminant reduces to
$$
(2(\alpha+\beta)(x+t)+A+B+2C)(2(\alpha-\beta)(x-t)+A+B-2C)\,.
$$
On $S_4$ the discriminant identically vanishes, so the singular set is all the plane and the corresponding tensors are not
characteristic tensors.
Outside of $S_4$,  if the discriminant is not constant (i.e., outside of $S_1\cap S_3$), then  $\Delta=0$ is a pair of orthogonal lines or a single line and the singular set is made of two opposite quadrants (the two webs corresponding to SC4) or of an half-plane (SC3). If $\Delta$ is a positive constant, the singular set is empty (SC1), while if it is negative all points are singular and the tensor is not characteristic ($(S_1 \cap S_3)-S_4$  not containing the plane $C=0$).
The classification given here can also be compared with that given in Table III
of \cite{McL}, where the type of any separable web in $\Min$ is characterized in terms
gamma and $I_\pm = \mathrm{sgn}\,(Z_\pm)$.  Note that in \cite{McL} (as in \cite{Kalnins}) the discrete
transformation $\widehat{R_2}$ is used, with the consequence that the number of distinct
types of separable webs is reduced from the ten described in the present paper
to nine.

\section{Conclusion}
We have classified Killing tensors of valence two in the Euclidean and Minkowski planes under the action of a group that preserves the type of the Killing web.  The method is based on a detailed analysis of the rank of the determining system of partial differential equations for the group invariants and depends crucially on the fact the generic rank of the system is six, which equals the dimension of the space of Killing two-tensors.  This result is dimensionally dependent.  It is thus unclear whether the method or a modification thereof can be extended to flat spaces of higher dimension or to spaces of non-zero constant curvature.  Nonetheless for the cases where the method is applicable it provides a very elegant algebraic classification for the type of the Killing web defined by a characteristic Killing tensor. 
This classification is equivalent to the classification of quadratic symmetric operators in the generators of the isometries of $\Min$, given in  \cite{Kalnins} and to the classification given in \cite{McL} in terms of Killing tensor
invariants, up to the exchange between space and time: since we do not allow a
change in signature of the metric, the coordinates of type SC4 (parabolic of type
I in \cite{Kalnins}) splits into the classes M7 and M8.  Our classification, not being
restricted to characteristic Killing tensors, extends the classification given
in \cite{McL} through the Invariant Theory of Killing Tensors.

\end{document}